\documentclass[11pt]{amsart}
\usepackage{amsmath,amsfonts,amsthm,amssymb,
amsthm,graphicx,mathtools,amscd}
\usepackage{hyperref}

\usepackage[mathscr]{eucal}
\usepackage{graphics}
\usepackage{color}
\usepackage{epsfig}

\usepackage[nocompress]{cite}

\DeclareSymbolFont{iwonaletters}{OML}{iwona}{m}{it}
\DeclareMathSymbol{\bdel}{\mathalpha}{iwonaletters}{"E}

\newtheorem{theorem}{Theorem}[section]
\newtheorem{lemma}[theorem]{Lemma}

\newtheorem{definition}[theorem]{Definition}

\newtheorem{proposition}[theorem]{Proposition}

\newtheorem{remark}[theorem]{Remark}

\newcommand{\beginsec}{
\setcounter{equation}{0}
}

\newcommand{\la}{\lambda}
\newcommand{\eps}{\varepsilon}

\newcommand{\al}{\alpha}

\newcommand{\kap}{\kappa}
\newcommand{\sig}{\sigma}
\newcommand{\del}{\delta}

\newcommand{\Gam}{\mathnormal{\Gamma}}
\newcommand{\Del}{\mathnormal{\Delta}}

\newcommand{\Om}{\mathnormal{\Omega}}

\newcommand{\N}{{\mathbb N}}

\newcommand{\R}{{\mathbb R}}

\newcommand{\Z}{{\mathbb Z}}

\newcommand{\EE}{{\mathbb E}}
\newcommand{\E}{{\mathbb E}}

\newcommand{\PP}{{\mathbb P}}

\newcommand{\calB}{{\mathcal B}}

\newcommand{\calF}{{\mathcal F}}

\newcommand{\calI}{{\mathcal I}}
\newcommand{\calJ}{{\mathcal J}}

\newcommand{\calM}{{\mathcal M}}
\newcommand{\calN}{{\mathcal N}}

\newcommand{\calP}{{\mathcal P}}
\newcommand{\calR}{{\mathcal R}}
\newcommand{\calS}{{\mathcal S}}

\newcommand{\calZ}{{\mathcal Z}}

\newcommand{\bX}{{\mathbf X}}

\newcommand{\frD}{\mathfrak{D}}

\newcommand{\frS}{\mathfrak{S}}

\newcommand{\frf}{\mathfrak{f}}

\renewcommand{\proof}{\noindent{\bf Proof.\ }}

\newcommand{\lan}{\langle}
\newcommand{\ran}{\rangle}

\newcommand{\w}{\wedge}

\newcommand{\To}{\Rightarrow}

\newcommand{\iy}{\infty}
\newcommand{\up}{\uparrow}
\newcommand{\loc}{{\rm loc}}
\newcommand{\cadlag}{c\`adl\`ag }
\newcommand{\osc}{{\rm osc}}

\newcommand{\noi}{\noindent}
\newcommand{\ds}{\displaystyle}

\newcommand{\rank}{{\rm rank}}

\newcommand{\ser}{{\rm ser}}
\newcommand{\dep}{{\rm dep}}

\newcommand{\stat}{{\rm stat}}
\newcommand{\mac}{{\rm mac}}
\newcommand{\jmp}{{\rm jmp}}
\newcommand{\cts}{{\rm cts}}

\newcommand{\Ll}{\mathbb{L}}

\begin{document}

\title[Invariance principle for load balancing]{
Invariance principle and McKean-Vlasov limit
for randomized load balancing in heavy traffic
}

\author{Rami Atar}
\address{Viterbi Faculty of Electrical and Computer Engineering
\\
Technion -- Israel Institute of Technology
} 
\email{rami@technion.ac.il}

\author{Gershon Wolansky}
\address{Department of Mathematics
\\
Technion -- Israel Institute of Technology
}
\email{gershonw@technion.ac.il}

\subjclass[2010]{68M20, 90B22, 60F17, 35K20}
\keywords{Randomized load balancing, power of choice,
low sampling rate, diffusion limit, hydrodynamic limit,
McKean-Vlasov limit, propagation of chaos,
parabolic initial boundary value problem, viscous scalar conservation law,
Daley-Miyazawa semimartingale representation.
}

\date{April 14, 2024}

\begin{abstract}
We consider a load balancing model where a Poisson
stream of jobs arrive at a system of many servers whose service time distribution
possesses a finite second moment.
A small fraction of arrivals pass through
the power-of-choice algorithm, which assigns a job to the shortest
among $\ell$, $\ell\ge2$, randomly chosen queues, and the remaining jobs
are assigned to queues chosen uniformly at random.
The system is analyzed at critical load in an asymptotic regime where
both the number of servers and the usual heavy traffic parameter
associated with individual queue lengths grow to infinity.
The first main result is a hydrodynamic limit stating that
the empirical measure of the diffusively normalized
queue lengths converges to a path in measure space
whose density is given by the unique solution of a parabolic PDE
with nonlocal coefficients.
The PDE has a stationary solution expressed explicitly,
which in particular provides a quantification of the balance achieved by the algorithm.
For fixed $\ell$, as the intensity of the load balancing stream
varies between its limits, this solution varies from exponential distribution
to a Dirac distribution, demonstrating that the result
covers the whole range from independence to state space collapse.

Further, two forms of an invariance principle are proved,
one under a rather general initial distribution and the other
under exchangeability of the initial conditions.
In the latter, the limit of individual normalized queue length is given by
a McKean-Vlasov SDE, and propagation of chaos holds.
The McKean-Vlasov limit is closely related to limit results
for Brownian particles on $\R_+$ interacting through their rank
(with a specific interaction). However, an entirely different set of tools is required,
as the collection of $n$ prelimit particles does not obey
a Markovian evolution on $\R_+^n$.
\end{abstract}

\maketitle

\section{Introduction}

This paper is concerned with a load balancing model
where a Poisson stream of arrivals faces a system of
$n$ servers working in parallel, with a general common
service time distribution having a finite second moment.
Following a setting proposed in \cite{banerjee2023load},
motivated by maintaining low communication overhead,
the stream undergoes thinning, where a small fraction of the arrivals
are routed via the {\it join the shortest of $\ell$ queues} (abbreviated JSQ($\ell$))
algorithm, while the other arrivals undergo {\it uniform routing}.
Here, JSQ($\ell$), also known as {\it power of choice}, is a randomized load balancing
algorithm that assigns a job to the shortest among $\ell$, $\ell\ge2$,
queues chosen uniformly at random, whereas uniform routing refers to
assigning a job to a queue chosen uniformly at random.

\begin{figure}
\begin{center}
\includegraphics[width=22em]{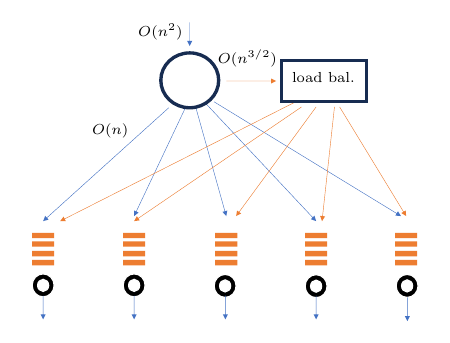}
\end{center}
\vspace{-1.6em}
\caption{\sl
Orders of magnitude of the various streams.
}
\label{fig2a}
\end{figure}

The focus is on the asymptotics as the number of servers, $n$,
becomes large, while individual queues fluctuate at the diffusion scale.
This is achieved by taking $n$ to also serve as the heavy traffic parameter
of individual queue lengths and critically loading the system. More precisely,
it is assumed that the per-server stream of arrivals resulting
from the uniform routing is asymptotic to $\la n+\hat\la n^{1/2}$,
whereas the processing rate of each server is asymptotic to $\la n+\hat\mu n^{1/2}$,
where $\la>0$, $\hat\la,\hat\mu\in\R$ are constants.
Moreover, the overall intensity of the load balancing
stream is taken asymptotic to $bn^{3/2}$, where $b>0$ is a constant
(see Fig.~\ref{fig2a}).
Note that the average load on a queue, defined as the rate of arrival of jobs per server
divided by a server's processing rate, is asymptotic to
$1+\rho n^{-1/2}$, where $\rho=\la^{-1}(\hat\la+b-\hat\mu)$.
In particular, $\rho<0$ is a stability condition.

Let $X^n_i(t)$, $i\in[n]:=\{1,\ldots,n\}$ denote the $i$-th queue length process
at the $n$-th system. Then the diffusion scaled queue length processes are given by
$\hat X^n_i(t)=n^{-1/2}X^n_i(t)$.
Let also $\bar\xi^n_t=n^{-1}\sum_{i\in[n]}\del_{\hat X^n_i(t)}$ denote
the empirical distribution of normalized queue lengths at time $t$,
where $\del_x$ is the Dirac measure at $x\in\R$.
In particular, $\bar\xi^n_t[0,x]$ is the fraction of queues $i$ for which $\hat X^n_i(t)\in[0,x]$.
As argued informally below, taking the load balancing stream intensity to scale
like $n^{3/2}$ is the only choice that gives interesting limits for $\bar\xi^n$.

The parameters $b$ and $\ell$
determine the volume of communication between the servers
and the dispatcher. They also
have dramatic influence on the degree to which
the system is load balanced. One of the goals of this work is
to quantify the level of achieved load balancing, which could be defined,
for example, by the expected empirical variance, namely
\begin{equation}\label{300}
\begin{split}
\sig^n(t)^2 :=&
\,\EE\Big[\frac{1}{n}\sum_{i\in[n]}\hat X^n_i(t)^2-\Big(\frac{1}{n}\sum_{i\in[n]}\hat X^n_i(t)\Big)^2\Big]
\\
=&\,\E\Big[\int_{\R_+}x^2\bar\xi^n_t(dx)-\Big(\int_{\R_+}x\bar\xi^n_t(dx)\Big)^2\Big].
\end{split}
\end{equation}
Whereas it is difficult to estimate $\sig^n(t)$
for each fixed $n$, suppose that the limit as $n\to\iy$ of $\bar\xi^n=\{\bar\xi^n_t\}$ exists
and is given by a deterministic path $\xi=\{\xi_t\}$ in the space of measures on $\R_+$.
Suppose, moreover, that $\xi_t$ has density $u(\cdot,t)$ for every $t$, i.e.,
$\xi_t(dx)=u(x,t)dx$.
This path could then be regarded as the macroscopic description of the model.
Thus one could instead look at the substitute
\begin{equation}\label{301}
\sig_\mac(t)^2:=\int_{\R_+}x^2u(x,t)dx-\Big(\int_{\R_+}xu(x,t)dx\Big)^2
\end{equation}
and obtain an approximation to $\sig^n(t)$ for large $n$.

Assuming $\bar\xi^n_0\to\xi_0$ in probability,
$\xi_0$ a Borel probability measure on $\R_+$,
our first main result (Theorem \ref{th2}) addresses this by providing a hydrodynamic limit,
which shows that $\bar\xi^n_t$ converges in probability,
as a measure-valued process, to a deterministic
measure-valued path given by $\xi_t(dx)=u(x,t)dx$, $t>0$. Here, the density $u$
is the unique classical solution of the initial-boundary value problem
\begin{equation}\label{14_}
\begin{cases}
u_t=-[(b\ell v^{\ell-1}-c_1)u]_x+au_{xx}
&(x,t)\in(0,\iy)^2\\
\ds
v(x,t)=\int_x^\iy u(y,t)dy
&(x,t)\in(0,\iy)^2\\
(c_1-b\ell)u(0,t)+au_x(0,t)=0
&t\in(0,\iy)\\
u(\cdot,t)dx\to\xi_0(dx)
&\text{weakly as $t\downarrow0$.}
\end{cases}
\end{equation}
Above, $c_1=-\la\rho+b$, and $a>0$ is a variance parameter.
The first line of \eqref{14_} contains a non-local drift coefficient that captures
the algorithm's instantaneous effect on the macroscopic distribution. Its third line
is a boundary condition expressing zero flux through the origin
(a Robin boundary condition; see e.g.\ \cite[p.\ 109]{pinchover2005introduction}), ensuring
that the total mass is preserved.
Our results also shows that $\sig^n(t)\to\sig_\mac(t)$ for each $t>0$ (Proposition \ref{prop2}).

When $\rho<0$, \eqref{14_} has a unique stationary solution given by
\begin{equation}\label{69}
\begin{split}
v_\stat(x)&=w(x)^{-1/(\ell-1)},
\qquad
u_\stat(x)=\frac{c_1}{a}(1-\al)e^{\frac{c_1}{a}(\ell-1)x}w(x)^{-\ell/(\ell-1)},
\\
w(x)&=(1-\al)e^{\frac{c_1}{a}(\ell-1)x}+\al,
\qquad
\al=\frac{b}{c_1}=\Big(1-\la\frac{\rho}{b}\Big)^{-1}\in(0,1),
\end{split}
\end{equation}
for $x\in\R_+$ (as stated below in Proposition \ref{prop1}).
While an explicit expression for the stationary distribution is interesting by its own right,
it also sheds light on the following aspect. If we fix $(\ell,\rho,\la,a)$
and let $b\to0$ and, respectively, $b\to\iy$, then the measure $u_\stat(x)dx$ converges
weakly to the exponential distribution and to $\del_0(dx)$
(details on this calculation are given in \S\ref{sec23};
graphs of $u_\stat$ for different values of $\ell$ and $b$
are shown in Fig.\ \ref{fig-var}).
These two limits indicate that our setting
covers the entire range from queues operating independently
to a {\it state space collapse}. The latter term is used in the heavy traffic
literature to describe a coordination between the queues to the degree
that the queue lengths, when normalized at the diffusion scale,
are asymptotically indistinguishable. In addition, this argument confirms the heuristic that
taking the load balancing intensity to scale like $n^{3/2}$ is the only choice for obtaining
nontrivial limits of $\bar\xi^n$. (At the same time, it suggests that
working with load balancing intensity of a higher order,
the behavior of $\{X^n_i/\eta(n)\}_{i\in[n]}$, for a suitable
choice of $1\ll\eta(n)\ll n^{1/2}$
may exhibit interesting behavior;
see \cite{gam16} for estimates on the order of magnitude of $X^n_i$ when all arrivals
are routed via JSQ($\ell$)).

Although we prove that \eqref{69} is the unique stationary solution
of the dynamics given by \eqref{14_}, let us emphasize that
in this paper we do not prove that every solution to \eqref{14_}
converges to $u_\stat$ as $t\to\iy$, nor that the interchange of the $t$ and $n$
limits is valid. These important aspects of the model will be the subject
of future work.

Our second goal is to provide an invariance principle.
This is proved under two different assumptions about initial conditions.
The first such result (Theorem \ref{th3}) assumes that, for a fixed $k$,
$(\hat X^n_i(0))_{i\in[k]}\To (X_i(0))_{i\in[k]}$, and states that
$(\hat X^n_i)_{i\in[k]}\To(X_i)_{i\in[k]}$ where the latter is the solution
to the SDE in $\R_+^k$
\begin{equation}
\label{100}
X_i(t)=X_i(0)+b_1t+b_0\int_0^tv(X_i(s),s)^{\ell-1}ds+\sig W_i(t)+L_i(t),
\qquad i\in[k],
\end{equation}
in which $L_i$ are boundary terms having continuous
nondecreasing sample paths, satisfying the condition $\int_{[0,\iy)}X_i(t)dL_i(t)=0$.
Here, $b_0=b\ell$, $b_1=-c_1$ and $\sig^2=2a$.
Alternatively (Theorem \ref{th4}),
if the initial queue lengths are assumed exchangeable,
then the limit in distribution of the $k$-tuple is given by $k$ independent
copies of
\begin{equation}
\label{101}
X(t)=X(0)+b_1t+b_0\int_0^tv(X(s),s)^{\ell-1}ds+\sig W(t)+L(t),
\end{equation}
with $\int_{[0,\iy)}X(t)dL(t)=0$.
In both cases, the diffusion limits depend on the service time distribution
only through its first two moments.

The second form of the invariance principle, where exchangeability
is assumed, combined with the characterization of the limiting
empirical distribution in terms of \eqref{14_}, shows that the pair $(X,\bar F)$,
describing the limiting stochastic dynamics of a rescaled queue length
and its law, satisfies
\begin{equation}\label{102}
\begin{split}
&X(t)=X(0)+\int_0^t\mathbf{b}(X(s),\bar F(\cdot,s))ds+\sig W(t)+L(t),
\\
&\bar F(x,t)=\PP(X(t)>x),
\end{split}
\end{equation}
where $\mathbf{b}(x,\bar F)=b_1+b_0\bar F(x)^{\ell-1}$.
Viewed this way, the result is closely related to the literature on interacting
diffusion models,
and in particular to the subject of diffusions interacting through their ranks,
where \eqref{102} is referred to as a McKean-Vlasov SDE.
For background on McKean-Vlasov limits and propagation of chaos we
refer to \cite{sznitman1991topics, gartner1988mckean}.

To explain the relation, consider a parametric regime
in which the heavy traffic parameter
is taken to its limit first. In this case, each normalized queue length
behaves as a Brownian particle on $\R_+$, with
an interaction among the particles caused by the load balancing algorithm.
The algorithm effects the dynamics by selecting a particle with probability
depending on its rank. In the original model, the selected queue
increases by one job. In the particle system this can be imitated
by imposing a positive rank-dependent drift.
It therefore comes as no surprise that the McKean-Vlasov SDE \eqref{102}
is a special case of the one that arises in the study of
rank-dependent diffusions. Indeed, existence and uniqueness of solutions to
\eqref{102}, as well as convergence results, follow from work on interacting
diffusions, specifically \cite[Sec.\ 2.4]{lacker2018strong}.

Clearly, the above explanation is only heuristic.
In the regime under consideration, the prelimit objects
are queue lengths, which do not follow a Markovian evolution
when considered as an $n$-tuple with state space $\R_+^n$.
As the basic assumptions of the interacting diffusion setting do not hold,
the treatment requires a completely different set of tools.

\subsection{Motivation and related work}

\ 

{\it Motivation: Invariance principle.}
Invariance principles for queueing models in heavy traffic
are diffusion limit results in which the limit processes depend on
model ingredients
such as service time and inter-arrival distributions only through their first two moments.
Their significance stems from the fact that they express robustness of performance to
underlying distributions.
As far as load balancing algorithms in many-server settings are concerned,
invariance principles have not been established before
and in fact, diffusion limit results have not been proved beyond
the case of exponential service times (for treatments of nonexponential service
in other asymptotic regimes, see below). This is
a rather restrictive assumption from a modeling viewpoint.

{\it Motivation: Non-perfect balancing, sparse messaging.}
A large body of work on load balancing
has been devoted to subcritically loaded systems, where on average queues
are short. Critically loaded systems have also been studied, where either state space
collapse or near optimal performance were proved; see references below.
The question of quantifying the degree of balancing achieved by
the algorithm, when balancing is far from
perfect, does not arise in such settings, and has not been addressed before
except under a fixed number of servers \cite{banerjee2023load}.
As already mentioned, this is one of the aspects of practical
significance that motivates this work.
In addition, the practical importance of load balancing
algorithms that operate under low communication volume
is widely acknowledged (see e.g.\ \cite{mendelson2022care})
and provides further motivation for our setting.

{\it Related asymptotic settings.}
The literature on load balancing in asymptotic regimes is vast;
see \cite{muk22sur} for a recent survey.
The present paper is not the first to study diffusion limits in
a regime in which the number of servers
scales like $n$ and individual queues are critically loaded,
with each queue length scaling like $n^{1/2}$. This precise parameterization
was considered before in \cite{bud-coh}, which studied
a multi-agent game rather than a load balancing model.

Closer to the present paper is \cite{gam16} which studies a system of $n$
exponential servers where all arrivals are routed via JSQ($\ell$),
under critical load. Rather than the empirical measure of rescaled queue lengths,
\cite{gam16} is concerned with a rescaling of the empirical measure
of the (unscaled) queue lengths, focusing of the short queues,
and identifying their dynamics under the limit.


The asymptotic regime where individual queues undergo heavy traffic
scaling and JSQ($\ell$) is applied on a small fraction of the arrival stream
was introduced in the aforementioned \cite{banerjee2023load}, in a setting
where the number of servers is fixed. This could model two different
scenarios, that are also relevant for our model.
One is where a single stream is split by the dispatcher.
Another is where each server has a dedicated stream
due, for example, to geographical or compatibility constraints,
and an additional stream is shared by all servers.
Under an exponential service assumption, it was proved in
\cite{banerjee2023load} that the collection of rescaled queue lengths converges
to a rank-based diffusion.

{\it Asymptotics of JSQ($\ell$).}
Perhaps the most well known load balancing algorithm is
the {\it join-the-shortest-queue} (JSQ), which routes jobs the shortest among all queues.
Its variation JSQ($\ell$), that is typically applied with $\ell$ much smaller that $n$,
significantly reduces the communication overhead while maintaining good performance.
It was introduced in \cite{vved96}, where under the
exponential service time assumption and subcritical load, the $n\to\iy$ limit
of the length of a typical queue in stationarity was shown to have
a doubly exponential
tail decay, a dramatic improvement over routing uniformly at random
in which case decay is exponential.
Related results appeared in \cite{mitz01}.
The result was extended in \cite{graham2000chaoticity} where
empirical measures were considered and propagation of chaos was
obtained. Functional central limit theorems (CLT) were obtained in \cite{graham2005functional} and
strong approximation results, including law of large numbers (LLN)
and CLT, in \cite{luczak2005strong}.
(These CLT results were not invariance principles and were not concerned
with the asymptotic regime studied in this paper).
Results on the mixing rate and the size of the maximal queue length were
obtained in \cite{luczak2006maximum}, and
aspects of stability and performance under
server heterogeneity were investigated in \cite{mukhopadhyay2015analysis}.
The paper \cite{muk18un} studied how
$\ell(n)$ should grow so that JSQ$(\ell(n))$ with
exponential servers would perform like JSQ, and thus achieve
asymptotic delay optimality.

A series of papers \cite{bramson2010randomized, bra12, bra13}
extended \cite{vved96} to several families of general service time distributions,
based on an approach that first establishes propagation of chaos
and then uses it to compute queue length distribution in equilibrium.
Another line of research treating JSQ($\ell$) with general service times is
\cite{aghajani2017pde, aghajani2019hydrodynamic}
and \cite{agarwal2020invariant}.
In \cite{aghajani2017pde, aghajani2019hydrodynamic},
the dynamical behavior was studied under general load and
the hydrodynamic limit
was shown to be given by an infinite system of coupled PDE.
A numerical method was developed to solve these PDE.
In \cite{agarwal2020invariant}, an infinite system of PDE was constructed
and shown to constitute the invariant state of the aforementioned system
of PDE under a subcriticality condition.
Our setting, where individual queues undergo heavy traffic
scaling, as well as our results, are quite different from both these
lines of work.

{\it Asymptotics of JSQ and other load balancing algorithms.}
JSQ is known to be delay optimal under exponential service
and asymptotically delay optimal in various limiting regimes;
see, for example, \cite{chen2012asymptotic}.
The diffusion limit in heavy traffic of JSQ in a many-server setting
was established in \cite{gam18}, where, assuming
exponential servers, a rescaled empirical measure
of queue lengths was shown to converge
to a process expressed in terms of a 2d diffusion.
Convergence of the steady state at the same scale
was proved in \cite{braverman2020steady}, and properties
of the diffusion process were investigated in
\cite{banerjee2019join, banerjee2020join}.
Under the regime considered in this line of work, the number of queues
that are of length $0$, $1$ and $2$ is of order $n^{1/2}$, $n$ and,
respectively, $n^{1/2}$, and only a negligible fraction exceeds length $2$.
Thus this regime captures quite a different behavior from what is described in this paper. The papers \cite{gupta2019load} and \cite{zhao2021many} considered
JSQ in other diffusive regimes, namely
the so called nondegenerate slowdown regime and, respectively,
the super-Halfin-Whitt regime.

Numerous load balancing algorithms besides JSQ and JSQ($\ell$)
have been proposed. Ones that emphasize sparse
communication, where the messaging rate often goes far below
that of JSQ($\ell$), include
\cite{van2020zero, ying2017power, sto17, atar2020persistent}.
Further asymptotic results on load balancing with nonexponential service times
include join-the-idle-queue \cite{sto17},
pull-based load distribution \cite{stolyar2015pull, stolyar2017pull},
zero waiting algorithms \cite{liu2022large},
and join-the-shortest-estimated-queue \cite{atar2021heavy}.

{\it Rank-based diffusions.}
McKean-Vlasov limits of diffusions interacting through their
ranks were addressed in \cite{shkolnikov2012large}.
Observing that in these models the dependence
of each particle's coefficients on the empirical law is discontinuous,
a situation not covered by the classical treatment such as
\cite{gartner1988mckean},
this paper gave convergence results
assuming that the coefficients are merely measurable.
The paper \cite{lacker2018strong} mentioned above addresses
a setting of which rank-based interaction is a special case,
and attains convergence results in a stronger topology than weak convergence.

\subsection{Notation}
Denote $\R_+=[0,\iy)$.
Let $\iota:\R_+\to\R_+$ denote the identity map.
In $\R^N$, denote the Euclidean norm by $\|\cdot\|$.
For $(\bX,d_\bX)$ a Polish space, let $C(\R_+,\bX)$ and $D(\R_+,\bX)$
denote the space of continuous and, respectively, \cadlag paths, endowed
with the topology of uniform convergence on compacts and, respectively,
the $J_1$ topology. Denote by $\calM_1$ the space of probability measures
on $\R_+$ equipped with the topology of weak convergence.
Denote by $C^\up$ the set of members of $C(\R_+,\R_+)$
that are nondecreasing and start at $0$.
For $\xi\in D(\R_+,\R^N)$, an interval $I\subset\R_+$, and $0\le \del\le T$, denote
\begin{align*}
\Del\xi(t)&=\xi(t)-\xi(t-) \qquad t>0, \qquad \Del\xi(0)=\xi(0),\\
\osc(\xi,I)&=\sup\{\|\xi(s)-\xi(t)\|:s,t\in I\},
\\
w_T(\xi,\del)&=\sup\{\|\xi(t)-\xi(s)\|:s,t\in[0,T],|s-t|\le\del\},
\\
\|\xi\|^*_T&=\sup\{\|\xi(t)\|:t\in[0,T]\},
\end{align*}
and by $|\xi|(t)$ the total variation of $\xi$ in $[0,t]$.

For $\frD\subset\R^N$, denote by $\frD^o$ and $\bar\frD$ its interior and,
respectively, closure. Denote by
$C(\frD)$ the set of continuous functions $f:\frD\to\R$ and by $C_0(\frD)$
the set of $f\in C(\frD)$ whose support is a compact subset of $\R^N$. For $k,l\in\N$ and
an open set $\frD\subset\R$ (resp., $\frD\subset\R\times\R_+$), denote by $C^k(\frD)$
(resp., $C^{k,l}(\frD)$) the set of functions $f:\frD\to\R$ possessing continuous
derivatives up to and including $k$ (resp., $(k,l)$).
For $\frD\subset\R$, $C^k(\frD)$ denotes the set of functions $f:\frD\to\R$
whose restriction to $\frD^o$ lies in $C^k(\frD^o)$, and whose derivatives up to
and including $k$ have continuous extensions to $\frD$. Define $C^{k,l}(\frD)$
analogously. For $\frD\in\R$ and $\frD\in\R\times\R_+$
denote $C^\iy(\frD)=\bigcap_{k\in\N}C^k(\frD)$ and
$C^\iy(\frD)=\bigcap_{k,l\in\N}C^{k,l}(\frD)$, resp. A subscript $0$ denotes compact
support, e.g.\ $C^k_0(\frD)=C^k(\frD)\cap C_0(\frD)$.
A subscript $b$ denotes boundedness, e.g.\ $C^k_b(\frD)$ are functions
in $C^k(\frD)$ whose derivatives of order $0\le l\le k$ are bounded.
For $f$ defined on (a subset of) $\R$, $f'$ denotes derivative,
and for $f$ defined on (a subset of) $\R\times\R_+$, $f_x$ and $f_t$ denote
spatial and temporal derivatives, resp.

For an interval $I\subset\R_+$ and a normed space $S$, denote by
$\Ll^p(I;S)$, $1\le p\le\iy$, the usual $\Ll^p$ space
defined in terms of the Lebesgue measure on $I$.
Denote by $\Ll^p_\loc(\R_+;S)$
the set of $f:\R_+\to S$ such that $f\in\Ll^p([0,T];S)$ for all finite $T$.
$\Ll^p(\R_+;\R)$ is abbreviated to $\Ll^p(\R_+)$.

For $f,g:\R_+\to\R$ and a measure $m$ on $\R_+$,
$\lan f,m\ran=\int_{\R_+}f(x)m(dx)$ and $\lan f,g\ran=\int_{\R_+}f(x)g(x)dx$.

If $V\in\R^n$ then $V_i$, $i\in[n]$ denote its components in the standard basis,
and vice versa: Given $V_i$, $i\in[n]$,
$V$ denotes the vector $(V_1,V_2,\ldots,V_n)$. Both these conventions hold
also for random variables $V_i$ and processes $V_i(\cdot)$.
Occasionally, with a slight abuse of standard terminology,
a sequence of random elements (random variables or processes) is referred to
as {\it tight} when their laws form a tight sequence of probability 
measures.
The term {\it with high probability} (w.h.p.) means ``holds, for each $n$,
on $\Om_n$, where $\lim_n\PP(\Om_n)=1$''.
Convergence in distribution is denoted by $\To$.
$c$ denotes a positive constant whose value may change from one expression
to another.

Throughout the paper, superscript $n$ attached to scalars, random variables
or processes, denotes dependence on the index $n$
rather than power.

\subsection{Paper organization}
The load balancing model is presented in \S\ref{sec21}. In \S\ref{sec22},
equation \eqref{111} is introduced, which is a PDE closely related to \eqref{14_},
and the main results are stated.
A discussion of the results appears in \S\ref{sec23}.
An outline of the proof is given in \S\ref{sec24}.
\S\ref{sec3} is devoted to developing several tools required for the proof.
The first tool, provided in \S\ref{sec31}, is a crucial PDE uniqueness result
for weak solutions of \eqref{111}.
Next, in \cite{miya19}, a semimartingale representation
for counting processes was introduced, that we find extremely useful.
\S\ref{sec32} presents this representation and, based on it, identifies
various martingales that are of importance to the model. These martingales
are key to the proof of convergence provided later in \S\ref{sec4}.
\S\ref{sec33} provides various estimates on the rescaled queueing processes.
Tightness of the empirical process $\bar\xi^n$ is shown in \S\ref{sec34}.
With this set of tools, the proof is then carried out in \S\ref{sec4},
which starts by writing down an equation for $\lan\phi,\bar\xi^n\ran$,
for $\phi$ a test function, namely equation \eqref{24+}.
Thanks to the uniqueness of solutions to \eqref{111}
and the tightness of $\bar\xi^n$, the hydrodynamic can be established
by showing that subsequential limits of $\bar\xi^n$ form weak solutions
to this PDE. This is achieved by relating limits of the various terms in
\eqref{24+} to weak solutions of \eqref{111}, performed in three steps
in \S\ref{sec41}, \S\ref{sec42} and \S\ref{sec43} (a detailed description of these
steps appears at the end of \S\ref{sec24}). Finally,
all the main results are proved, based on the above, in \S\ref{sec44}.

\section{Load balancing in heavy traffic}\label{sec2}
\beginsec

\subsection{The load balancing model}\label{sec21}

\subsubsection{Arrivals, queue lengths and basic relations}

In the model there are $n$ servers and a queue in front of each.
There is a dedicated stream of arrivals into each queue
and an additional stream of arrivals, called the {\it load balancing stream}
(LBS), that go through
the JSQ$(\ell,n)$ algorithm. These $n+1$ arrival streams are modeled
as mutually independent Poisson processes.
Clearly, this could be recast as a single Poisson stream,
out of which $n+1$ thinned streams are created
by means of random selection.

In what follows, the queueing systems will be indexed by $n\in\N$,
the number of servers.
The processes $X^n_i$, $E^n_i$, $D^n_i$ and $T^n_i$ represent
the $i$-th queue length process, dedicated arrival process,
departure process and busyness process, respectively.
Denote by $A^n_0$ the LBS arrival process,
and by $A^n_i$ the process counting LBS arrivals routed to server $i$.
For each $i$, $E^n_i$ is a Poisson process of parameter
$\la^n$, and $A^n_0$ is Poisson of parameter $\la^n_0$,
all having right-continuous sample paths. The $n+1$ Poisson processes
are mutually independent, for each $n$. We have
\begin{equation}\label{01}
X^n_i(t)=X^n_i(0-)+E^n_i(t)+A^n_i(t)-D^n_i(t),
\qquad i\in[n],\ t\in\R_+.
\end{equation}
Work conservation is assumed, hence
\begin{equation}\label{03}
T^n_i(t)=\int_0^t1_{\{X^n_i(s)>0\}}ds.
\end{equation}

\subsubsection{The load balancing algorithm}

Upon each LBS arrival,
$\ell$ out of the $n$ queues, chosen uniformly at random, are sampled.
The arrival is routed to the queue that is shortest among the $\ell$;
if there are ties, the queue with the smaller index is preferred.
Given $x\in\R^n$, let
\begin{equation}\label{17-}
\rank(i;x)=\#\{j:x_j<x_i\}+\#\{j\le i:x_j=x_i\}.
\end{equation}
Then the probability that a LBS arrival
is routed to the queue whose rank, defined by \eqref{17-},
is $r$ is given by
\begin{equation}
\label{19}
p_{n,r}=\frac{{n-r\choose \ell-1}}{{n\choose \ell}},
\qquad r\in[n],
\end{equation}
with ${k\choose j}=0$ when $j>k$. Note that
\begin{equation}\label{17}
p_{n,r}=0 \text{ for } r\ge n-\ell+2,
\qquad
\max_{r\in[n]} p_{n,r}=p_{n,1}=\frac{\ell}{n}.
\end{equation}
Thus the randomization mechanism can equivalently be achieved
by letting, for each $n$, $\theta^n_k$, $k\in\N$ be IID random variables
with $\PP(\theta^n_1=r)=p_{n,r}$, $r\in[n]$, and
routing the $k$-th LBS arrival to the queue
whose rank is $\theta^n_k$. That is,
\begin{equation}\label{02}
A^n_i(t)=\int_{[0,t]}1_{\{\calR^n_i(s-)=\theta^n_{A^n_0(s)}\}}dA^n_0(s),
\qquad
\calR^n_i(t)=\rank(i;X^n(t)),
\qquad i\in[n],\ t\in\R_+.
\end{equation}
Of course, this randomization scheme is not used in practice as it
requires keeping track of the queue lengths of the entire system.
However, it is convenient for carrying out the analysis.
\begin{remark}\label{rem6}
Although we work with the expression \eqref{19} that corresponds to sampling without
replacement, all our results are valid for sampling with replacement as well.
This is stated and explained in Remark \ref{rem6+}.
\end{remark}

\subsubsection{The initial condition}
We allow quite a general initial condition,
where residual times of jobs already being processed
at time $0$ may be dependent and have unspecified
distributions. Denote the (random) set of queues that at time $0$ contain
no jobs and, respectively, at least one job, by
$\calN^n=\{i\in[n]:X^n_i(0)=0\}$ and $\calP^n=\{i\in[n]:X^n_i(0)>0\}$.
For $i\in\calP^n$, let $Z^n_i(0)$ denote the initial residual time
of the head-of-the-line job in queue $i$.
For $i\in\calN^n$ we add fictitious jobs having
zero processing time. To this end, rather than specifying
$X^n(0)$ as the initial queue length, $X^n(0-)$ is specified;
and for each $i\in\calN^n$,
the queue length is set to $X^n_i(0-)=1$ and the residual processing time
to $Z^n_i(0)=0$. Obviously, this results in $X^n_i(0)=0$.
This convention allows us to greatly simplify notation when we later construct
counting processes for service and departure.
The initial condition is thus a tuple
\[
\calI^n=(\{X^n_i(0-),Z^n_i(0),i\in[n]\},\calN^n,\calP^n),
\]
where $(\calN^n,\calP^n)$ partitions $[n]$, and
\[
\begin{split}
X^n_i(0-)&=1,\ Z^n_i(0)=0,\ i\in\calN^n,
\\
X^n_i(0-)&\ge1,\ Z^n_i(0)>0,\ i\in\calP^n.
\end{split}
\]

\subsubsection{Service times}

Let $\Phi_\ser$ be a Borel probability measure on $[0,\iy)$
with mean 1 and standard deviation $\sig_\ser\in(0,\iy)$, such that
$\Phi_\ser(\{0\})=0$. Let $\Phi^n_\ser$ be defined as
a scaled version of $\Phi_\ser$ uniquely specified via
$\Phi^n_\ser[0,x]=\Phi_\ser[0,\mu^nx]$, $x\in\R_+$.
Here, $\mu^n>0$ is the service rate in the $n$-th system.
For $k\ge1$, let $Z^n_i(k)$ denote the service time
of the $k$-th job to be served by server $i$ after the head-of-the-line
job at time $0-$ there (for $i\in\calN^n$ this means the $k$-th
job after the fictitious one). For every $i$, $\calZ^n_i=(Z^n_i(k),k\ge1)$
is an IID sequence with common distribution $\Phi^n_\ser$.

Next, the potential service process $S^n_i$, evaluated at $t$,
gives the number of jobs completed by server $i$ by the time it has
worked $t$ units of time. It is given, with $\sum_0^{-1}=0$, by
\[
S^n_i(t)=
\max\Big\{k\in\Z_+:\sum_{j=0}^{k-1}Z^n_i(j)\le t\Big\}, \qquad t\ge0.
\]
The departure processes are given by $D^n_i(t)=S^n_i(T^n_i(t))$.
This is the number of jobs completed by time $t$ by server $i$.
Note that the first departure counted by $D^n_i$
is the one initially processed if $i\in\calP^n$
and the fictitious job if $i\in\calN^n$.

\subsubsection{Dependence structure}

For each $n$, the $2n+3$ stochastic elements
\begin{equation}\label{095}
\text{$E^n_i$, $i\in[n]$, $\calZ^n_i$, $i\in[n]$, $\calI^n$,
$A^n_0$, and $\{\theta^n_k\}$
are mutually independent.}
\end{equation}

\subsubsection{Scaling and critical load condition}
The arrival and service rates are assumed to satisfy the following.
There are constants $\la>0$ and $\hat\la\in\R$ such that
\begin{equation}\label{04}
\hat\la^n:=n^{-1/2}(\la^n-n\la)\to\hat\la \text{ as } n\to\iy,
\end{equation}
a constant $b>0$ such that
\begin{equation}\label{05}
\hat\la^n_0:=n^{-3/2}\la^n_0\to b \text{ as } n\to\iy,
\end{equation}
and constants $\mu>0$ and $\hat\mu\in\R$ such that
\begin{equation}\label{06}
\hat\mu^n:=n^{-1/2}(\mu^n-n\mu)\to\hat\mu \text{ as } n\to\iy.
\end{equation}
The critical load condition is assumed, namely
\begin{equation}\label{09}
\la=\mu.
\end{equation}
Some further notation used throughout is
\[
\hat b^n_1=\hat\la^n-\hat\mu^n,\qquad b_1=\hat\la-\hat\mu,
\quad
c_1=-b_1,
\quad
b_0=b\ell,
\quad
\sig^2=\la(1+\sig_\ser^2),
\quad
a=\frac{\sig^2}{2}.
\]
Note that the average load on a queue,
defined as the rate of arrival of work per server,
$\la^n+n^{-1}\la^n_0$, divided by
a server’s processing rate, $\mu^n$, that is,
\[
\frac{n\la+n^{1/2}\hat\la+n^{1/2}b+o(n^{1/2})}{n\mu+n^{1/2}\hat\mu+o(n^{1/2})},
\]
is asymptotic to
\[
1+n^{-1/2}\rho, \qquad \rho:=\frac{\hat\la+b-\hat\mu}{\la}.
\]
Hence the {\it load parameter} $\rho$ is given by
$\rho=(b_1+b)/\la$.

Let rescaled versions of the queue length processes
and cumulative idle time processes be defined by
\begin{equation}\label{07}
\hat X^n_i(t)=n^{-1/2}X^n_i(t),
\qquad
\hat L^n_i(t)=n^{-1/2}\mu^n(t-T^n_i(t)),
\end{equation}
and denote
$\bar\xi^n_t=n^{-1}\sum_{i\in[n]}\del_{\hat X^n_i(t)}$.
Further, assume that $\bar\xi^n_{0-}\to\xi_0$
(equivalently $\bar\xi^n_0\to\xi_0$)
in probability, as $n\to\iy$, where $\xi_0$ is a deterministic
Borel probability measure on $\R_+$.
For the rescaled residual times $Z^n_i(0)$, assume
that for every $\eps>0$,
\begin{equation}\label{40}
\lim_{n\to\iy}\max_{i\in[n]}\PP(\tilde Z^n_i(0)>\eps)=0 \quad\text{ where}\quad
\ \tilde Z^n_i(0):=\tilde\mu^nZ^n_i(0),\ \tilde\mu^n:=n^{-1/2}\mu^n,
\end{equation}
and
\begin{equation}\label{62}
\sup_n\max_{i\in[n]}\EE[\tilde Z^n_i(0)^2]<\iy.
\end{equation}
Because $\tilde\mu^n$ scales like $n^{1/2}$, \eqref{40} and \eqref{62}
impose the condition that $Z^n_i(0)$ scale like $o(n^{-1/2})$.
This is mild compared to the condition on $Z^n_i(k)$, $k\ge1$,
which have been assumed to scale like $n^{-1}$.
Finally, assume that
\begin{equation}\label{53}
\sup_n\max_{i\in[n]}\EE[\hat X^n_i(0-)^2]<\iy.
\end{equation}
Note that, by Fatou's lemma, this imposes a condition on $\xi_0$,
namely $\xi_0$ necessarily satisfies $\int x^2\xi_0(dx)<\iy$.

All assumptions made thus far are in force throughout the paper.

\subsection{Main results}\label{sec22}

First we address well-posedness of the relevant PDE, starting
with classical solutions to \eqref{14_}.

\begin{theorem}\label{th1}
Within the class of functions $u\in C^{2,1}((0,\iy)\times\R_+;\R_+)$
satisfying, for every $T<\iy$, $\sup_{t\in(0,T]}\int_{\R_+} xu(x,t)dx<\iy$,
there exists a unique solution to equation \eqref{14_}.
Moreover, for each $t>0$, $u(\cdot,t)$ is a probability density.
\end{theorem}

A function $u_\stat\in C^2(\R_+)$ is said to be a {\it stationary solution
associated with \eqref{14_}}, if it is a probability density
possessing a finite second moment,
and, setting $\xi_0(dx)=u_\stat(x)dx$, the solution to \eqref{14_}
is given by $u(x,t)=u_\stat(x)$ for all $x,t$.

\begin{proposition}\label{prop1}
Assume $\rho<0$.
Then there exists a unique stationary solution to \eqref{14_}.
It is given by \eqref{69}.
\end{proposition}

Theorem \ref{th1} is a consequence of a result stated next,
that gives uniqueness of weak solutions to a class of
viscous scalar conservation laws, and provides
one of the main tools used in this paper. This is
a parabolic equation of the type
\begin{equation}\label{111}
\begin{cases}
v_t=\left(\frf(v)\right)_x+av_{xx}, & (x,t)\in\R_+^2,\\
v(0,t)=1, & t> 0,\\
v(\cdot,0)=v_0.
\end{cases}
\end{equation}
 A key point is that this result does not
require any regularity assumption in the $x$ variable on the class of solutions.
The definition of a weak solution is as follows.

\begin{definition}\label{222} 
A function 	$v\in\mathbb{L}_\loc^\infty(\R_+; \mathbb{L}^\infty(\R_+))\cap \mathbb{L}_\loc^\infty(\R_+; \mathbb{L}^1(\R_+))$ is a weak solution of (\ref{111}) if for
any $t\in (0,\iy)$ and any $\phi\in C^\infty_0(\R_+)$
satisfying $\phi(0)=0$,
\begin{align}\label{80}
\langle v(\cdot,t), \phi\rangle - \langle v_0, \phi\rangle
= -\int_0^t\langle \frf(v(\cdot,s)), \phi'\rangle  ds
+ a\int_0^t\langle v(\cdot,s), \phi''\rangle ds +a\phi'(0)t.
\end{align}
\end{definition}
The choice of $\frf$ that will be of interest is $\frf(z)=c_1z-bz^\ell$.
In particular, \eqref{14_} is related \eqref{111} when the latter takes the special form
\begin{equation}\label{111+}
\begin{cases}
v_t=(c_1v-bv^\ell)_x+av_{xx}, & (x,t)\in\R_+^2,\\
v(0,t)=1, & t> 0,\\
v(\cdot,0)=\xi_0(\cdot,\iy).
\end{cases}
\end{equation}
The existence and uniqueness of a classical solution of \eqref{111}
is well known, and there is a vast literature on weak solutions
in the $W^{1,p}_\loc$ sense. However, we found no uniqueness
result of a weak solution as given in Definition \ref{222} under the mere
assumption
$v\in\mathbb{L}_\loc^\infty(\R_+; \mathbb{L}^\infty(\R_+))\cap \mathbb{L}_\loc^\infty(\R_+; \mathbb{L}^1(\R_+))$.

\begin{theorem}\label{th1+}
Assume $\frf\in C^\infty(\R)$ and $v_0\in \mathbb{L}^1(\R_+)\cap\mathbb{L}^\infty(\R_+)$. 
Then there exists a unique weak solution to (\ref{111})
(in particular, to \eqref{111+}) in the sense
of Definition \ref{222}. This solution is in $C^\iy(\R_+\times(0,\iy))$.
If $v$ denotes the weak solution to \eqref{111+} then
$u=-v_x$ is the classical solution of \eqref{14_}.
\end{theorem}

The following is our first main result, showing that \eqref{14_} provides
a macroscopic description of the model, specifically its hydrodynamic limit.

\begin{theorem}\label{th2}
Let $\xi_0$ be extended to a trajectory $\xi=\{\xi_t,t\in\R_+\}$ in $\calM_1$
by setting
\[
\xi_t(dx)=u(x,t)dx, \qquad t>0,
\]
where $u$ is the unique solution to \eqref{14_}. Then $\xi\in C(\R_+,\calM_1)$,
and one has $\bar\xi^n\to\xi$ in probability in $D(\R_+,\calM_1)$,
as $n\to\iy$.
\end{theorem}

Under a slightly stronger moment condition, the above result implies
\begin{proposition}\label{prop2}
Assume that, for some $\eps>0$,
$\Phi_\ser$ possesses a finite $2+\eps$ moment and that
$\sup_n\max_{i\in[n]}\{\E[\tilde Z^n_i(0)^{2+\eps}]\vee\E[\hat X^n_i(0-)^{2+\eps}]\}<\iy$ (compare with \eqref{62} and \eqref{53}).
Let $\sig_\mac$ from \eqref{301} be defined with $u$ of \eqref{14_}.
Then $\sig^n(t)\to\sig_\mac(t)$ for every $t>0$.
\end{proposition}

Our second main result is the conjunction of following two invariance principles,
obtained under two different assumptions on the initial conditions.

\begin{theorem}\label{th3}
Fix $k\in\N$ and assume that $(\hat X^n_i(0))_{i\in[k]}\To(X_i(0))_{i\in[k]}$.
Let $v$ be the weak solution to \eqref{111+}.
Then $(\hat X^n_i,\hat L^n_i)_{i\in[k]}\To(X_i,L_i)_{i\in[k]}$
in $(D(\R_+,\R_+)\times C(\R_+,\R_+))^k$, where the latter tuple
is the unique in law solution to the system \eqref{100},
driven by a $k$-dimensional standard BM $(W_i)_{i\in[n]}$
independent of $(X_i(0))_{i\in[k]}$.
\end{theorem}

\begin{theorem}\label{th4}
Assume that, for each $n$, the random variables
$X^n_i(0-)$, $i\in[n]$, are exchangeable
(necessitating that the limit law of each $X^n_i(0-)$ is $\xi_0$).
Let $v$ be the weak solution to \eqref{111+}.
Then for $k\in\N$,
$(\hat X^n_i,\hat L^n_i)_{i\in[k]}\To(X_i,L_i)_{i\in[k]}$
in $(D(\R_+,\R_+)\times C(\R_+,\R_+))^k$, where the latter tuple
is given by $k$ independent copies of the unique in law solution
$(X,L)$ to \eqref{101}
driven by a standard BM $W$ independent of $X(0)$, and the latter is distributed
according to $\xi_0$.
\end{theorem}
In view of Theorem \ref{th2}, the SDE in Theorem \ref{th4}
is a McKean-Vlasov SDE, as it can be seen to take the form
\eqref{102}.
Note that we do not assume that the entire initial data is exchangeable,
nor that the tie breaking rule is symmetric w.r.t.\ exchanging $i$'s.
The result implies that these issues have negligible effect upon taking the limit.

\subsection{Discussion}\label{sec23}

\begin{figure}
\begin{center}
\includegraphics[width=20em]{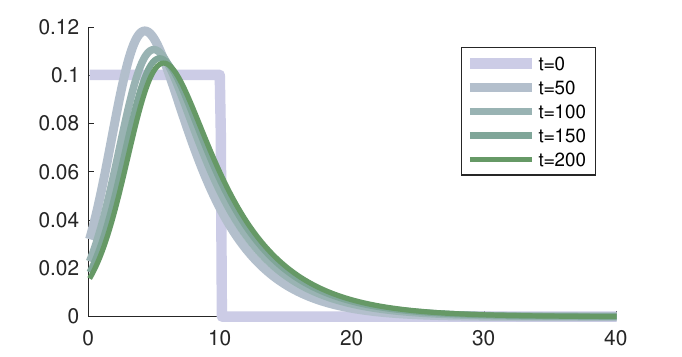}
\includegraphics[width=20em]{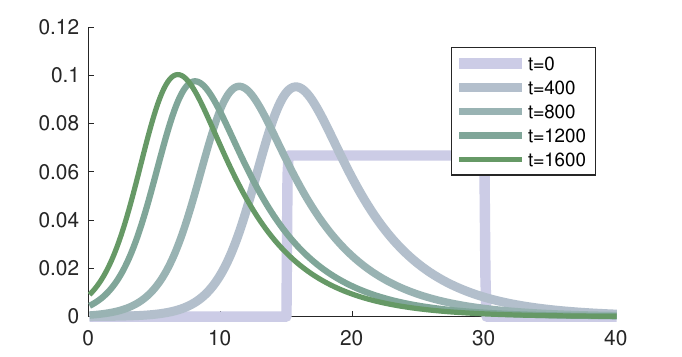}
\end{center}

\caption{\sl
Solution of \eqref{14_} for initial condition $\rm unif[0,10]$
(left) and $\rm unif[15,30]$ (right) at different time instances,
with $\rho=-.01$, $c_1=.21$, $b=.2$,
$\ell=4$, $a=1$.
}
\label{fig-ex}
\begin{center}
\includegraphics[width=20em,height=12em]{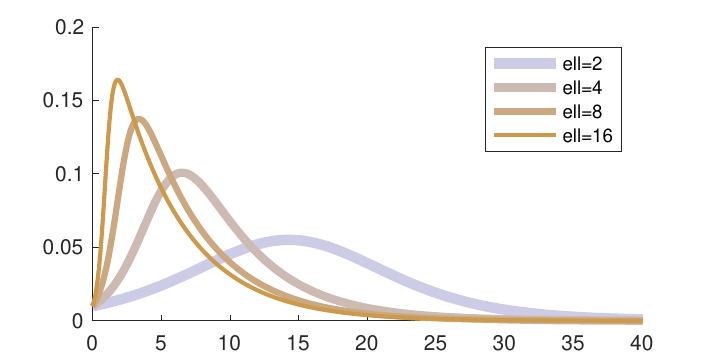}
\includegraphics[width=20em,height=12em]{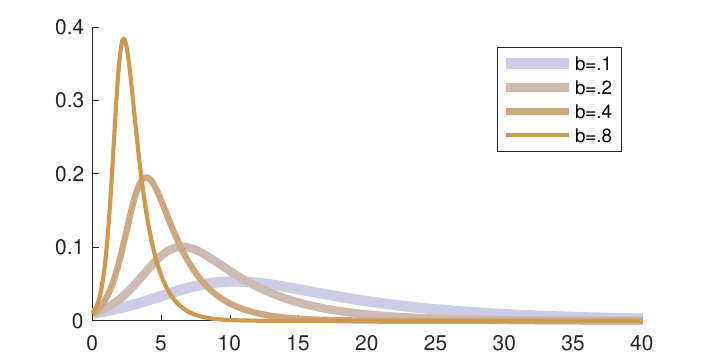}
\end{center}
\caption{\sl
Stationary solution for different values of $\ell$,
with $b=.2$ (left), and for different values of
$b$, with $\ell=4$ (right). In all cases, $\rho=-.01, a=1,\la=1$.
}
\label{fig-var}
\begin{center}
\includegraphics[width=20em, height=11em]{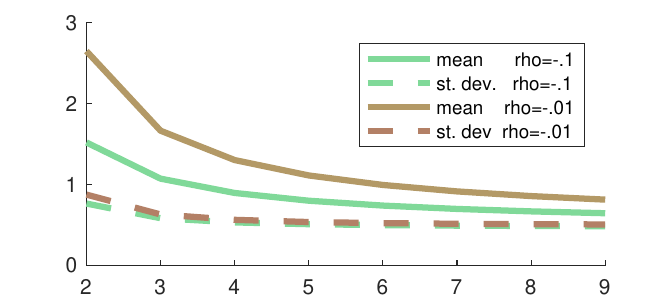}
\includegraphics[width=20em, height=11em]{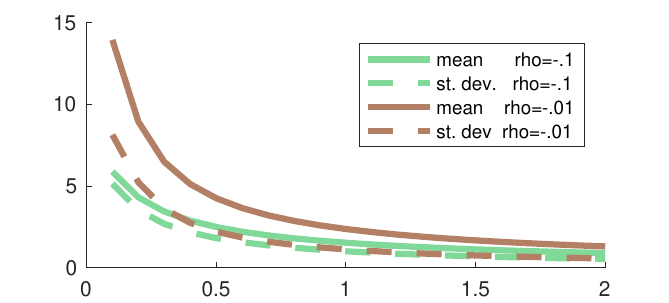}
\end{center}
\caption{\sl
Stationary mean and standard deviation for $2\le\ell\le9$,
with $b=2$, $a=1$ (left)
and for $0.1\le b\le2$, with $\ell=4$, $a=1$ (right).
}
\label{fig-mn-stdev}
\end{figure}

The hydrodynamic limit result establishes that the dynamics of the model
are given at the macroscopic scale by the PDE \eqref{14_}. As argued in the introduction,
one can use this to express basic performance criteria in terms of the solution $u$.
In particular, denoting the macroscopic mean by
\[
m_\mac(t):=\int_0^\iy xu(x,t)dx,
\]
a first order approximation to the mean queue length is given by $m_\mac(t)n^{1/2}$.
Accordingly, the mean delay is, to first order, given by $\la^{-1}m_\mac(t)n^{-1/2}$. Similarly, 
the macroscopic standard deviation, $\sig_\mac(t)$ defined in \eqref{301}, may be taken as
an {\it index of balance}.
As a second order parabolic equation in one spatial variable,
there exist standard tools for numerically solving \eqref{14_}
and evaluating these performance criteria.
Some examples of solutions are plotted in Fig.\ \ref{fig-ex}.

Whereas the macroscopic dynamics can be solved
numerically, the macroscopic equilibrium state has an explicit formula,
namely \eqref{69}.
As already mentioned, it is beyond the scope of this paper
to show rigorously that all solutions $u$ converge
to $u_\stat$ as $t\to\iy$, and further, that
the invariant distribution of the stochastic model's dynamics
converges to $u_\stat$ as $n\to\iy$.
However, even without establishing these results, $u_\stat$ constitutes a legitimate
solution to the PDE. In particular,
a combination of Theorem \ref{th2} and Proposition \ref{prop1}
implies that $\bar\xi^n\to\xi$ in probability holds provided that $\xi_0(dx)=u_\stat(x)dx$.
Here, $\xi_t=\xi_0$ for all $t$.
Fig.\ \ref{fig-var} shows
graphs of $u_\stat$ for different values of $\ell$ and $b$,
based on \eqref{69}.

We can also use $u_\stat$ in place of $u(\cdot,t)$
in the above macroscopic performance criteria, and define
analogously $m_\mac$ and $\sig_\mac$ which are now
time independent quantities.
Fig.\ \ref{fig-mn-stdev} shows their dependence on $\ell$ and $b$.

Let us provide details on the statement made in the introduction
regarding the limits, as $b\to0$ and $b\to\iy$, of the stationary distribution.
Fix $(\ell,\rho,\la,a)$, with $\rho<0$. According to the formula in \eqref{69},
as $b\to0$, one has $\al\to0$ and $c_1\to-\la\rho=\la|\rho|$.
Hence $w(x)\to e^{\frac{\la|\rho|}{a}(\ell-1)x}$ pointwise, and consequently,
\[
\lim_{b\to0}v_\stat(x)=e^{-\frac{\la|\rho|}{a}x},\qquad x\ge0.
\]
As $b\to\iy$, one has $c_1\to\iy$ and, although $\al\to1$, it is obvious
that $w(x)\to\iy$ for every $x>0$. This gives
\[
\lim_{b\to\iy}v_\stat(x)=1_{\{x=0\}}, \qquad x\ge0.
\]
In both cases, the limit holds at every continuity point of the limiting CDF,
proving that $v_\stat(dx)=u_\stat(x)dx\to \frac{\la|\rho|}{a}e^{-\frac{\la|\rho|}{a}x}dx$
weakly as $b\to0$, and
$u_\stat(x)dx\to\del_0(dx)$ weakly as $b\to\iy$.
The former is an exponential distribution, which agrees with the limiting distribution
of a single server queue in heavy traffic. This makes perfect sense since in this
limit the queues operate independently. The latter is a distribution that has
all its mass at $0$, showing that, by choosing $b$ large,
the hydrodynamic limit can be made
arbitrarily close to a state space collapse, where the diffusion-scale processes
$\hat X^n_i$ are asymptotically indistinguishable.

The nature of the drift in equation \eqref{14_}
is akin to the dynamics of a mass distribution subjected to a gravitational field
\cite{wol95},
where the force acting on an infinitesimal mass depends on the amount of mass
that piles above it.
Here it is the mass accumulated {\it below} it,
specifically, a nonlinear function of $v(x,t)$,
that determines the force.

\subsection{Outline of the proof}\label{sec24}

The approach used to proving the hydrodynamic limit result is
based on PDE uniqueness. Thus the first main tool required
is the uniqueness of weak solutions to \eqref{111}.
The crux of the argument is as follows.  If $v_1, v_2$ are two weak solutions
then one can show that $w=v_1-v_2$ satisfies the equality
$$
\langle w(\cdot,t), \psi(\cdot,t)\rangle = \int_0^t \Big\langle w(\cdot,s),\Big[ \psi_s(\cdot,s)- \Big[\frac{\frf(v_1)-\frf(v_2)}{v_1-v_2}\Big]\psi_x+\psi_{xx}\Big] 
\Big\rangle  ds
$$
for any $t>0$ and test function $\psi$ satisfying $\psi\in C_0^\infty(\mathbb{R}_+\times[0,t])$,   $\psi(0,s)=0$ for $s\in [0,t]$.  One can then find an $\mathbb{L}^\infty$ approximation of the solution  to the linear, backward equation
$$ \psi_s(\cdot,s)- \Big[\frac{f(v_1)-f(v_2)}{v_1-v_2}\Big]\psi_x+\psi_{xx}=0$$
 on $\mathbb{R}_+\times [0,t]$ satisfying $\psi(\cdot,t)=\phi(\cdot)$,
given any smooth $\phi$. As a result,
$\langle w(\cdot,t), \phi\rangle =0$ for any $t$ and any such $\phi$,
hence $w\equiv 0$. 

Next,  $\bar\xi^n$ is shown to form a tight sequence, and then the next main step
is to prove that limit points satisfy \eqref{111} in weak sense. Achieving this goal
relies on two key elements. One is a semimartingale decomposition for
point processes, introduced in \cite{miya19}. This decomposition
is particularly convenient in our setting, where the number of point processes involved,
given by the departure processes, grows to infinity. The other is an estimate showing
$C$-tightness of $\hat X^n_i$ uniformly in $i\in[n]$.
With this toolbox we can then show the following. Let $\tilde\phi$
be a test function as in Definition \ref{222} and $\phi$ its antiderivative.
Then every limit point $\xi$ of $\bar\xi^n$ satisfies
\begin{align}\label{70}
\lan \phi,\xi_t\ran &= \lan \phi,\xi_0\ran+\int_0^t\lan b_1\phi'+a\phi'',\xi_s\ran ds
+\frac{b_0}{\ell}\int_0^t\int_{\R_+}\phi'(x)\frS(\xi_s[x,\iy),\xi_s(x,\iy))\xi_s(dx)ds,
\end{align}
where
\begin{equation}\label{51}
\frS(a,b)=a^{\ell-1}+a^{\ell-2}b+\cdots+ab^{\ell-2}+b^{\ell-1},
\qquad a,b\in\R.
\end{equation}
To see the relation to equation \eqref{80} required by Definition \ref{222},
let $v(x,t)=\xi_t(x,\iy)$.
Then $\lan\phi,\xi_t\ran=\int_{\R_+}(v(x,t)-1)\tilde\phi(x)dx$.
Suppose that $\xi_t$ has no atoms for every $t>0$.
Then, owing to the fact that $\frS(a,a)=\ell a^{\ell-1}$,
\eqref{70} reduces to \eqref{80} when
$\frf(z)=c_1z-bz^\ell$, which is the case of interest.
However, an a priori proof of the atomless property is not required, as a calculation via integration by parts
shows that \eqref{70} and \eqref{80} are equivalent even in presence of atoms.
As a result, the existence of weak limit of $\bar\xi^n$ satisfying the PDE
follows, completing the proof of the hydrodynamic limit.

Equation \eqref{111} is well known to have $C^\iy$ solutions,
and as a consequence of the above result the relation of the hydrodynamic limit
to the PDE \eqref{14_} for the density follows.

Finally, the two invariance principles follow by combining the $C$-tightness
of individual rescaled queue lengths with the hydrodynamic limit result,
where the atomless property, that by now has been proved,
gives a control over the interaction term.

We now provide some more details on the structure of \S\S \ref{sec3}--\ref{sec4}.
In \S\ref{sec31} we prove uniqueness of weak solutions to \eqref{111}.
The aforementioned semimartingale representations
are provided in \S\ref{sec32}. In \S\ref{sec33}
uniform estimates on the rescaled queueing processes are derived,
including $C$-tightness
and second moment. Tightness of $\bar\xi^n$ is shown in \S\ref{sec34}.

In \S\ref{sec41} equation \eqref{24+} for $\lan\phi,\bar\xi^n\ran$ is written and
it is shown that its linear term converges to the linear term in \eqref{70}.
In \S\ref{sec42}, the interaction term is shown to converge to that in \eqref{70}.
The reduction, alluded to above, of \eqref{70} to \eqref{80}, is shown in \S\ref{sec43}.
In \S\ref{sec44}, the convergence of $\bar\xi^n$ to a solution of \eqref{80}
is deduced from the above steps, and, based on that,
the proof of all main results is then provided.

\section{Preliminaries}\label{sec3}
\beginsec

\subsection{PDE uniqueness}\label{sec31}

Here we prove the uniqueness part of Theorem \ref{th1+},
stated as Lemma~\ref{lem9} below.
In preparation for proving this result, we provide an extension of
the class of test functions allowed in Definition \ref{222},
which is relatively standard but is given for completeness.

\begin{lemma}\label{lem8}
If $v$ is a weak solution of \eqref{111} (in the sense of Definition \ref{222})  then 
for any $t\in (0,\iy)$ and any $\psi\in C_0^\infty\left(\R_+\times[0,t]\right)$
satisfying $\psi(0,s)=0$, $0\leq s \le t$, one has
\begin{align}
\notag
\label{333}
\langle v(\cdot,t), \psi(\cdot,t)\rangle- \langle v_0, \psi(\cdot,0)\rangle
&= \int_0^t  \langle v(\cdot,s), \psi_s(\cdot,s)\rangle ds
-\int_0^t\langle \frf(v(\cdot,s)), \psi_x\rangle  ds
\\
&\qquad
+ a\int_0^t\langle v(\cdot,s), \psi_{xx}\rangle ds
+a\int_0^t\psi_x(0,s)ds.
\end{align}
\end{lemma}

\proof
Observe by \eqref{80} that for any function $\phi$ satisfying Definition \ref{222}
and for any $0\le\tau_1<\tau_2\leq t$, 
\begin{equation}
\label{eq-a}
\langle v(\cdot,\tau_2), \phi\rangle = \langle v(\cdot, \tau_1), \phi\rangle
-\int_{\tau_1}^{\tau_2}\langle \frf(v(\cdot,s)), \phi'\rangle  ds
+ a\int_{\tau_1}^{\tau_2}\left( \langle v(\cdot,s), \phi''\rangle +\phi'(0)\right)ds.
\end{equation}
Since $v$ and $\frf(v)$ are in $\Ll^\infty(\R_+)$ uniformly in $t$,
it follows, in particular, that $t\mapsto \langle v(\cdot,t), \phi\rangle$ is continuous
and, moreover,
\begin{equation}
\label{eq-b}
|\langle v(\cdot,\tau_2)-v(\cdot, \tau_1), \phi\rangle| \leq |\tau_1-\tau_2| K
\end{equation}
where $K=K(\phi,v)$.

Let now $\psi\in C_0^\infty(\R_+\times[0,t])$ be a function
that satisfies the conditions of the lemma.
Let $N\in \mathbb{N}$. Let $t_k:= tk/N$ where $k=0,\ldots, N$. Define $\phi_k(x)=\psi(t_k, x)$.  Then, from (\ref{eq-a})
$$\langle v(\cdot,t_{k+1}), \phi_{k}\rangle - \langle v(\cdot, t_k), \phi_k\rangle =-\int_{t_k}^{t_{k+1}}\langle \frf(v(\cdot,s)), \phi'_k\rangle  ds+ a\int_{t_k}^{t_{k+1}}\left(\langle v(\cdot,s), \phi''_k\rangle +  \phi'_k(0)\right)ds.
$$ 
Summing over $k=0, \ldots, N-1$ and using $\phi_N=\psi(\cdot, t)$, $\phi_0=\psi(\cdot,0)$ we get
\begin{multline}\label{444}\langle v(\cdot,t), \psi(\cdot,t(1-1/N))\rangle - \langle v_0, \psi(\cdot,0)\rangle -\sum_{k=1}^{N-1}\langle v(\cdot,t_{k+1}), \phi_{k+1}-\phi_k\rangle  \\ 
 =\sum_{k=0}^{N-1}\left[-\int_{t_k}^{t_{k+1}}\langle \frf(v(\cdot,s)), \phi'_k\rangle  ds+ a\int_{t_k}^{t_{k+1}}\left(\langle v(\cdot,s), \phi''_k\rangle +\phi'_k(0)\right)ds \right].
\end{multline}
Evidently, $\langle v(\cdot,t), \psi(\cdot,t(1-1/N))\rangle\rightarrow \langle v(\cdot,t), \psi(\cdot,t)\rangle$ 
as $N\rightarrow\infty$. 
Next, by (\ref{eq-b})
\begin{align*}
\langle v(\cdot,t_{k+1}), \phi_{k+1}-\phi_k\rangle
&= \int_{t_k}^{t_{k+1}} \langle v(\cdot,s), \psi_s(\cdot,s)\rangle ds+ \int_{t_k}^{t_{k+1}} \langle v(\cdot,t_{k+1})-v(\cdot, s), \psi_s(\cdot,s)\rangle ds
\\
&= \int_{t_k}^{t_{k+1}} \langle v(\cdot,s), \psi_s(\cdot,s)\rangle ds +O(N^{-2}).
\end{align*}
Thus
$$ \sum_{k=1}^{N-1}\langle v(\cdot,t_{k+1}), \phi_{k+1}-\phi_k\rangle \rightarrow \int_0^t \langle v(\cdot,s), \psi_s(\cdot,s)\rangle ds$$
as $N\rightarrow\infty$. 
Likewise, the right side of (\ref{444}) converges to
$$ -\int_0^{t}\langle \frf(v(\cdot,s)), \psi_{x}\rangle  ds+ a\int_{0}^{t}\left(\langle v(\cdot,s), \psi_{xx}\rangle+\psi_x(0,s)\right)\rangle ds. $$
This proves the lemma.
\qed

\begin{lemma}\label{lem9}
Assume $\frf\in C^\infty(\R)$ and $v_0\in \mathbb{L}^1(\R_+)\cap\mathbb{L}^\infty(\R_+)$. 
Then there is at most one weak solution to (\ref{111}) in the sense
of Definition \ref{222}.
\end{lemma}

\proof
As can be seen by performing a change of variables $x\mapsto a^{-1/2}x$,
we may and will assume w.l.o.g.\ that $a=1$.
Suppose $v_1, v_2$ are two solutions. The goal is to show that
$v_1=v_2$.
Let $w=v_1-v_2$ and note that $w\in\Ll^\iy_\loc(\R_+,\Ll^\iy(\R_+))
\cap\Ll^\iy_\loc(\R_+,\Ll^1(\R_+))$.
Let
$$J(x,t):=\begin{cases}
\ds
\frac{\frf(v_1(x,t))-\frf(v_2(x,t))}{v_1(x,t)-v_2(x,t)}  & \text{if} \ \ w(x,t)\not= 0  \\ 
\frf'(v_1(x,t) )  & \text{if} \ \ w(x,t)=0.
\end{cases}
$$ 
Then for any test function $\psi$ satisfying the conditions of Lemma \ref{lem8},
\begin{equation}\label{223}
\langle w(\cdot,t), \psi(\cdot,t)\rangle = \int_0^t  \langle w(\cdot,s),\left[  \psi_s(\cdot,s)- J(\cdot,s)\psi_x+\psi_{xx}\right] 
\rangle  ds.
\end{equation}
Fix $T$.
By the assumptions on $\frf$ and the definition of a solution,
$J\in \mathbb{L}^\infty([0, T]; \mathbb{L}^\infty(\R_+))$.
Moreover, for any $x_0>0$, $|\frf'(x)-\frf'(0)|\le cx$ holds provided that
$|x|<x_0$, with $c=c(x_0)$. Hence
$J-\frf'(0)\in \mathbb{L}^\infty([0, T];\mathbb{L}^1(\R_+))$. 
Let now $J_N\in C^\infty(\R_+\times[0,T])$
be a sequence that is bounded uniformly in $\Ll^\iy(\R_+\times[0,T])$
and satisfies $\|J_N-J\|_1=O(N^{-1})$, where we denote $\|\cdot\|_1=
\|\cdot\|_{\Ll^1(\R_+\times[0,T])}$.
Given $t\in(0,T]$,
let $\tilde\psi^N$ be the classical solution of
the backward linear problem on the time interval $[0,t]$:
 \begin{equation}\label{005}  \tilde\psi^N_s(\cdot,s)- J_N(\cdot,s)\tilde\psi^N_x+\tilde\psi^N_{xx}=0, \ \ 
\tilde\psi^N(\cdot,t)=\phi, \ \ \tilde\psi^N(0,s)=0, \qquad 0\le s<t,
\end{equation}
where $\phi\in C^\infty_0(\R_+)$.
Since $\tilde\psi^N$ and all its derivatives decay to zero as $x\rightarrow\infty$,
uniformly in $s\in[0,t]$, we may replace them by $\psi^N$ which 
satisfy  the  conditions of Lemma \ref{lem8} and
$$\|\psi^N-\tilde\psi^N\|_\infty+ \|\psi^N_x-\tilde\psi^N_x\|_\infty+\|\psi^N_{xx}-
\tilde\psi^N_{xx}\|_\infty +\|\psi^N_t-\tilde\psi^N_t\|_\iy = O(N^{-1})$$ 
uniformly on $[0,t]$. It follows from \eqref{223} that 
\begin{align*}
\langle w(\cdot,t), \phi\rangle &= \int_0^t  \langle w(\cdot,s),[\tilde\psi^N_s(\cdot,s)- J(\cdot,s)\tilde\psi^N_x+\tilde\psi^N_{xx}] 
\rangle  ds
\\
&=  \int_0^t  \langle w(\cdot,s),\left[  J_N(\cdot,s)-J\right]\tilde \psi^N_x
\rangle  ds
\\
&=  \int_0^t  \langle w(\cdot,s),\left[  J_N(\cdot,s)-J\right]\psi^N_x
\rangle  ds +O(N^{-1})
\end{align*}
for any such $\phi$. 
Suppose
\begin{equation}\label{224}
\sup_N\sup_{s\in[0,t]}\|\tilde\psi^N_x(\cdot,s)\|_\infty<\iy
\end{equation}
if $t$ is small enough.
Since $w\in\Ll^\iy([0,T];\Ll^\iy(\R_+))$ and
$\|J_N-J\|_1=O(N^{-1})$, this implies $\langle w(\cdot,t), \phi\rangle =O(N^{-1})$
for any $\phi\in C^\iy_0(\R_+)$ and any $N$, hence $w(\cdot,t)=0$. As a consequence,
$w(\cdot,t)=0$ for all $t\in[0,t_0]$, some $t_0>0$.
Moreover, $t_0$ does not depend on the initial condition $v_0$.
Thus, iterating the argument shows that $w=0$.

It thus suffices to show \eqref{224}. To this end,
denote $M := \sup_N\|J_N\|_\infty<\infty$. 
Let $\psi^0$ be the solution of 
$$  \psi^0_s(\cdot,s)+\psi^0_{xx}=0, \ \ 
 \psi^0(\cdot,t)=\phi, \ \ \psi^0(0,s)=0,  \ \  0\leq s\leq t. $$
Let $m_N(y,\tau):= J_N(y, \tau)\tilde\psi^N_y(y, \tau)$.
Since, as mentioned earlier, $\|\tilde\psi^N_y\|_\iy<\iy$ for all $N$,
one also has $\|m_N\|_\iy<\iy$ for all $N$.
Moreover, $\tilde\psi^N=\psi_0 + \hat{\psi}^N$ where $\hat\psi^N$ is the solution of 
$$ \hat\psi^N_s+\hat\psi^N_{xx} = m_N, \ \ \ \hat\psi^N(\cdot,t)=0, \ \hat\psi^N(0,s)=0, \qquad 0\leq s\leq t. $$
The solution $\hat\psi^N$ is given by Duhamel's principle: 
$$\hat\psi^N(x,s)=\frac{1}{2\sqrt{\pi}}\int_s^t d\tau  (t-\tau)^{-1/2}\left[\int_0^\infty m_N(y, t-\tau) e^{-\frac{(x-y)^2}{4(t-\tau)} }dy
	-\int_0^\infty m_N(y, t-\tau) e^{-\frac{(x+y)^2}{4(t-\tau)} }dy\right].$$
	 Then $\hat\psi^N_x(x,s)=-A_N(x,s)+B_N(x,s)$, where
\begin{align*}
A_N(x,s)&=\frac{1}{4\sqrt{\pi}}\int_s^t d\tau (t-\tau)^{-3/2}\int_0^\infty m_N(y, t-\tau) (x-y)e^{-\frac{(x-y)^2}{4(t-\tau)} }dy,
\\
B_N(x,s)&=\frac{1}{4\sqrt{\pi}}\int_s^t d\tau (t-\tau)^{-3/2}\int_0^\infty m_N(y, t-\tau) (x+y)e^{-\frac{(x+y)^2}{4(t-\tau)} }dy.
\end{align*}
Both $|A_N(x,s)|$ and $|B_N(x,s)|$ are bounded by
\[
\frac{\|m_N\|_\iy}{4\sqrt{\pi}}\int_s^td\tau(t-\tau)^{-3/2}
\int_{-\iy}^\iy|x+y|e^{-\frac{(x+y)^2}{4(t-\tau)}}dy.
\]
Changing variables of integration,
$$
|A_N(x,s)|\vee|B_N(x,s)|\le
\frac{ \|m_N\|_\infty}{2\sqrt{\pi}}\int_s^t d\tau  (t-\tau)^{-1/2}
\int_{-\iy}^\infty |z|e^{-\frac{z^2}{2 }}dz \leq C  \|m_N\|_\infty t^{1/2}, $$
for a constant $C$.
Recall that $\|m_N\|_\infty\leq M \|\psi^N_x\|_\infty$
and $\tilde\psi^N=\psi^0+\hat\psi^N$. 
Thus
$$\|\tilde\psi^N_x(\cdot,s)\|_\infty
\leq \|\psi^0_x\|_\infty+ 2C\|m_N\|_\infty t^{1/2}
\leq \|\psi^0_x\|_\infty+ 2CM \|\psi^N_x(\cdot,s)\|_\infty t^{1/2}\ , \qquad 0<s\le t.
$$
This shows $\|\tilde\psi^N_x(\cdot ,s)\|_\infty\le 2\|\psi^0_x\|_\iy$
for $0\leq s\leq (2CM)^{-2}/2$ and completes the proof that $v_1=v_2$.
\qed

\subsection{Martingale toolbox}\label{sec32}

A semimartingale decomposition of
counting processes was introduced in \cite{miya19}, which deviates
from the classical Doob-Meyer decomposition and is convenient
for our purposes. For the renewal process $S^n_i$, this decomposition
is constructed as follows.
Denote by
\[
R^n_i(t)=\inf\{s>0:S^n_i(t+s)>S^n_i(t)\}, \qquad t\ge0,
\]
the residual time to the next counting instant
(note that it is right-continuous, hence at a time of counting
it already shows the time until the next counting).
One has
\begin{equation}\label{36}
t+R^n_i(t)=\sum_{k=0}^{S^n_i(t)}Z^n_i(k).
\end{equation}
Hence with
\[
M^{\ser,n}_i(t)=\sum_{k=1}^{S^n_i(t)}\zeta^n_i(k),
\qquad \zeta^n_i(k)=1-\mu^nZ^n_i(k),\qquad k\ge1,
\]
and $\sum_1^0=0$, one has the
{\it Daley-Miyazawa semimartingale representation},
\begin{equation}\label{37}
S^n_i(t)=\mu^n(t-Z^n_i(0)+R^n_i(t))+ M^{\ser,n}_i(t),\qquad t\ge0.
\end{equation}
Denoting $\calF^{\ser,n}_i(t)=\sig\{S^n_i(s),R^n_i(s),s\le t\}$, the
first term on the right is adapted to the filtration $\{\calF^{\ser,n}_i(t)\}$
while $M^{\ser,n}_i$ a martingale on it. A form of this decomposition
that will be useful here is based on a filtration $\calF^n_t$ that, for each $t$,
contains all relevant information about the system by time $t$.
It will be obtained once a time change transformation of \eqref{37}
is performed, to represent $D^n_i(t)=S^n_i(T^n_i(t))$.

To this end, some further notation is required.
In addition to the delayed renewal processes
$S^n_i$ it will also be useful to introduce the corresponding non-delayed
renewal processes
\[
S^{0,n}_i(t)=
\max\Big\{k\in\Z_+:\sum_{j=1}^{k-1}Z^n_i(j)\le t\Big\}, \qquad t\ge0.
\]
Note that $S^{0,n}_i(0)=1$ and that these processes
are IID (unlike $S^n_i$).
Additional rescaled processes are denoted as follows
\begin{equation}\label{07+}
\hat E^n_i(t)=n^{-1/2}(E^n_i(t)-\la^nt),
\quad
\hat S^n_i(t)=n^{-1/2}(S^n_i(t)-\mu^nt),
\quad
\hat S^{0,n}_i(t)=n^{-1/2}(S^{0,n}_i(t)-\mu^nt),
\end{equation}
\begin{equation}\label{08}
\hat A^n_0(t)=n^{-1/2}(A^n_0(t)-\la^n_0t),
\quad
\hat A^n_i(t)=n^{-1/2}A^n_i(t).
\end{equation}
Next, let
\[
\tilde R^n_i(t)=n^{-1/2}\mu^nR^n_i(t)
=\tilde\mu^n R^n_i(t),
\qquad
\hat M^{\ser,n}_i(t)=n^{-1/2}M^{\ser,n}_i(t),
\]
and
\begin{equation}\label{35}
M^{\dep,n}_i(t)=\sum_{k=1}^{D^n_i(t)}\zeta^n_i(k),
\qquad
\hat M^{\dep,n}_i(t)=n^{-1/2}M^{\dep,n}_i(t).
\end{equation}
Then, by \eqref{37},
\begin{equation}\label{54}
D^n_i(t)=\mu^n(T^n_i(t)-Z^n_i(0)+R^n_i(T^n_i(t)))+M^{\dep,n}_i(t),
\end{equation}
and
\begin{equation}\label{34}
\hat S^n_i(T^n_i(t))=-\tilde Z^n_i(0)+\tilde R^n_i(T^n_i(t))
+\hat M^{\dep,n}_i(t).
\end{equation}
The tuple
\[
\calS^n(t)=(E^n_i(t),A^n_i(t),D^n_i(t),X^n_i(t),T^n_i(t),
R^n_i(T^n_i(t)),\,i\in[n],\, A^n_0(t),\theta^n_{A^n_0(t)})
\]
is referred to as the {\it state of the system at time $t$}. Denote
the corresponding filtration by
\[
\calF^n_t=\sig\{\calI^n,\calS^n(s),s\in[0,t]\}.
\]
Note that at times $t$ when server $i$ is active,
$R^n_i(T^n_i(t))$ is the residual time till the completion
of service of the job being processed by that server, whereas
at times when the server is idle, $R^n_i(T^n_i(t))$ gives
the service duration of the job that will be processed next by this server.

\begin{lemma}\label{lem1}
i.
The processes $\hat M^{\dep,n}_i$, $\hat E^n_i$, and $\hat A^n_0$
are $\{\calF^n_t\}$-martingales, with optional quadratic variations given by
\[
[\hat M^{\dep,n}_i](t)=n^{-1}\sum_{k=1}^{D^n_i(t)}\zeta^n_i(k)^2,
\qquad
[\hat E^n_i](t)=n^{-1}E^n_i(t),
\qquad
[\hat A^n_0](t)=n^{-1}A^n_0(t),
\]
and one has
\begin{equation}\label{64}
\EE\{[\hat M^{\dep,n}_i](t)\}=n^{-1}\sig_\ser^2\EE\{D^n_i(t)\}<\iy.
\end{equation}

\noi
ii. For distinct $i,j\in[n]$,
\begin{equation}\label{38-}
\E\{[\hat M^{\dep,n}_i,\hat M^{\dep,n}_j](t)\}=0.
\end{equation}
iii. The process
\[
\hat M^{A,n}_i(t)=\hat A^n_i(t)-\hat C^{A,n}_i(t),
\quad \text{ where } \quad
\hat C^{A,n}_i(t)=\la^n_0 n^{-1/2}\int_0^tp_{n,\calR^n_i(s)}ds,
\]
is an $\{\calF^n_t\}$-martingale, whose optional quadratic variation is
$[\hat M^{A,n}_i](t)=n^{-1}A^n_i(t)$.
\end{lemma}

\proof
i. For adaptedness of $\hat M^{\dep,n}_i$ it suffices to prove
that $Z^n_i(k)1_{k\le D^n_i(t)}\in\calF^n_t$ for all $k$.
This is shown as follows.
By \eqref{36},
$T^n_i(s)+R^n_i(T^n_i(s))=\sum_{k=0}^{D^n_i(s)}Z^n_i(k)$.
Hence $\{Z^n_i(k),k\le D^n_i(t)\}$ can all be recovered from
the tuple $\{T^n_i(s), R^n_i(T^n_i(s)), D^n_i(s)\}$, as $s$ varies
between $0$ and $t$.
Since the latter is $\{\calF^n_t\}$-adapted, this proves the claim.

Next it is shown that $\hat M^{\dep,n}_i(t)\in L_1(d\PP)$.
Note first that as a renewal process, $S^n_i(t)$ has finite expectation
for every $t$. Since $T^n_i(t)\le t$ this gives $\E[D^n_i(t)]<\iy$.
Let
\[
t^n_i(k)=\inf\{t\ge0:D^n_i(t)\ge k\},
\qquad
k=1,2,\ldots.
\]
These are clearly stopping times on $\{\calF^n_t\}$. Hence
\begin{equation}\label{38}
t^n_i(k)\in\calF^n_{t^n_i(k)-}, \qquad k\ge1,
\end{equation}
where we recall that for a stopping time $\tau$,
\[
\calF^n_{\tau-}=\calF^n_0\vee\sig\{A\cap\{\tau<t\}:A\in\calF^n_t,t\ge0\}
\]
(see \cite[I.1.11 and I.1.14]{jacshi}).
The state of the system up to $t^n_i(k)-$, namely
$\{\calS^n(t),t<t^n_i(k)\}$, can be recovered from the tuple
$\calI^n$, $(E^n_i(t),t\in\R_+,i\in[n])$, $(A^n_0(t), t\in\R_+)$,
$(\theta^n_j, j\in\N)$,
$(Z^n_\ell(j), j\in\N)$, $\ell\in[n]\setminus\{i\}$ and
$(Z^n_i(j), j<k)$, as follows by the construction of the model.
By our assumptions, $Z^n_i(k)$ is independent of this tuple.
As a result, it is independent of $\calF^n_{t^n_i(k)-}$.
It follows from \eqref{38} that
\[
\{D^n_i(t)\le k\}=\{t^n_i(k)\ge t\}\in\calF^n_{t^n_i(k)-}.
\]
The structure that we have just proved, where
\begin{equation}\label{63}
Z^n_i(k')\in\calF^n_{t^n_i(k)-}, k'< k, \text{ whereas $Z^n_i(k)$
is independent of $\calF^n_{t^n_i(k)-}$},
\end{equation}
the fact that $D^n_i(t)$ is a stopping time on the discrete parameter filtration
$\{\calF^n_{t^n_i(k)-},k\in\N\}$,
$\E[D^n_i(t)]<\iy$, along with the fact that $Z^n_i(k)$ are IID
with $\E[Z^n_i(k)]<\iy$, allows us to use Wald's identity, showing
\[
\E\Big[\sum_{k=1}^{D^n_i(t)}Z^n_i(k)\Big]
=\EE\{D^n_i(t)\}\EE\{Z^n_i(1)\}=\EE\{D^n_i(t)\}(\mu^n)^{-1}<\iy.
\]
Since $|\zeta^n_i(k)|\le 1+\mu^nZ^n_i(k)$, this shows
that $\EE\{\hat M^{\dep,n}_i(t)\}<\iy$.

To show the martingale property, note that by the independence stated in \eqref{63},
we have
\begin{equation}\label{39}
\EE[\zeta^n_i(k)|\calF^n_{t^n_i(k)-}]=0,\qquad k\ge 1.
\end{equation}
Arguing now along the lines of the proof of \cite[Lemma 2.1]{miya19},
\begin{align*}
\hat M^{\dep,n}_i(t)&=n^{-1/2}\sum_{k=1}^{D^n_i(t)}\zeta^n_i(k)
=n^{-1/2}\sum_{k=1}^\iy\zeta^n_i(k)1_{\{t^n_i(k)\le t\}}.
\end{align*}
Hence for $s<t$,
\begin{align*}
\EE[\hat M^{\dep,n}_i(t)|\calF^n_s]-\hat M^{\dep,n}_i(s)
&=
n^{-1/2}\sum_{k=1}^\iy\EE[\zeta^n_i(k)1_{\{s<t^n_i(k)\le t\}}|\calF^n_s]
\\
&=
n^{-1/2}\sum_{k=1}^\iy\EE[\EE[\zeta^n_i(k)
1_{\{s<t^n_i(k)\le t\}}|\calF^n_{t^n_i(k)-}]\,|\calF^n_s]
\\
&=0,
\end{align*}
where we used \eqref{38} and \eqref{39}.

For the martingale property of $\hat E^n_i$ one only needs to show that $E^n_i(t)-E^n_i(s)$
is independent of $\calF^n_s$ when $s<t$.
Again, this follows from the fact that all the processes comprising $\calS^n_u$,
$u\le s$, can be recovered from the tuple
$\calI^n$, $(E^n_i(u),u\in[0,s],i\in[n])$, $(A^n_0(u),u\in[0,s])$,
$(\theta^n_k, k\in\N)$, $(Z^n_i(k), k\in\N,i\in[n])$;
but the increment $E^n_i(t)-E^n_i(s)$ is independent of this tuple.

A similar proof holds for $\hat A^n_0$.

The expressions for the quadratic variation are straightforward.

To show \eqref{64} we can use Wald's identity as before,
now with the IID sequence $\zeta^n_i(k)^2$, which now gives
\[
\EE\Big[\sum_{k=1}^{D^n_i(t)}\zeta^n_i(k)^2\Big]
=\EE\{D^n_i(t)\}\EE\{\zeta^n_i(1)^2\},
\]
and \eqref{64} follows.

ii.
Similar to the argument following \eqref{38},
one can recover the state of the system up to $t^n_{ijkl}-$, where
\[
t^n_{ijkl}=\min(t^n_i(k),t^n_j(l)),
\]
namely $\{\calS^n(t),t<t^n_{ijkl}\}$, from
$\calI^n$, $(E^n_i(t),t\in\R_+,i\in[n])$, $(A^n_0(t), t\in\R_+)$,
$(\theta^n_j, j\in\N)$,
$(Z^n_{i'}(k), k\in\N)$, $i'\in[n]\setminus\{i,j\}$,
$(Z^n_i(k'), k'<k)$ and $(Z^n_j(l'),l'<l)$.
However, by our assumptions, the pair
$(Z^n_i(k),Z^n_j(l))$ is independent of this tuple,
hence it is independent of $\calF^n_{t^n_{ijkl}-}$.
Since $Z^n_i(k)$ and $Z^n_j(l)$ are mutually independent, this gives
\begin{equation}\label{399}
\EE[\zeta^n_i(k)\zeta^n_j(l)|\calF^n_{t^n_{ijkl}-}]=0,
\qquad k,l\ge 1.
\end{equation}
We have
\begin{align*}
[\hat M^{\dep,n}_i,\hat M^{\dep,n}_j](t)
&=n^{-1}\sum_{k=1}^{D^n_i(t)}\zeta^n_i(k)\sum_{l=1}^{D^n_j(t)}
1_{\{t^n_j(l)=t^n_i(k)\}}\zeta^n_j(l)
\\
&=n^{-1}\sum_{k=1}^\iy\sum_{l=1}^\iy\zeta^n_i(k)\zeta^n_j(l)
1_{\{t^n_j(l)=t^n_i(k)\le t\}}.
\end{align*}
In view of \eqref{38},
$1_{\{t^n_j(l)=t^n_i(k)\le t\}}\in\calF^n_{t^n_{ijkl}-}$.
Hence by \eqref{399},
\[
\EE[\zeta^n_i(k)\zeta^n_j(l)1_{\{t^n_j(l)=t^n_i(k)\le t\}}|\calF^n_{t^n_{ijkl}-}]=0,
\qquad
k,l\ge1,
\]
and \eqref{38-} follows.

iii. Clearly $A^n_i$, defined in \eqref{02} is adapted and $A^n_i(t)$ is
integrable for all $t$. Hence the same is true for $\hat M^n_i$. Next, let
\[
s^n(k)=\inf\{t\ge0:A^n_0(t)\ge k\},\qquad k=1,2,\ldots.
\]
As in (i), these are stopping times and $s^n(k)\in\calF^n_{s^n(k)-}$.
To show the martingale property, we can write, using
$\int_0^tp_{n,\calR^n_i(s-)}ds=\int_0^tp_{n,\calR^n_i(s)}ds$,
\begin{align*}
n^{1/2}\hat M^{A,n}_i(t)&=\int_{[0,t]}1_{\{\calR^n_i(s-)=\theta^n_{A^n_0(s)}\}}dA^n_0(s)
-\la^n_0\int_0^tp_{n,\calR^n_i(s)}ds
\\
&=\int_{[0,t]}(1_{\{\calR^n_i(s-)=\theta^n_{A^n_0(s)}\}}
-p_{n,\calR^n_i(s-)})dA^n_0(s)
+\int_0^tp_{n,\calR^n_i(s-)}(dA^n_0(s)-\la^n_0ds)
\\
&=: M^n_{i,1}(t)+M^n_{i,2}(t).
\end{align*}
For $M^n_{i,1}$, write
\[
A^n_i(t)=\sum_{k=1}^{A^n_0(t)}
1_{\{\calR^n_i(s^n(k)-)=\theta^n_k\}}.
\]
An argument similar to the one given before shows that
$\theta^n_k$ is independent of $\calF^n_{s^n(k)-}$.
Therefore, for $0\le s<t$, we have
\begin{align*}
\E[A^n_i(t)|\calF^n_s]-A^n_i(s)
&=\sum_{k=1}^\iy\E[1_{\{\calR^n_i(s^n(k)-)=\theta^n_k\}}
1_{\{s<s^n(k)\le t\}}|\calF^n_s]
\\
&=
\sum_{k=1}^\iy\E[\E[1_{\{\calR^n_i(s^n(k)-)=\theta^n_k\}}
1_{\{s<s^n(k)\le t\}}|\calF^n_{s^n(k)-}]|\calF^n_s]
\\
&=
\sum_{k=1}^\iy\E[p_{n,\calR^n_i(s^n(k)-)}
1_{\{s<s^n(k)\le t\}}|\calF^n_s]
\\
&=
\E[C^n_i(t)|\calF^n_s]-C^n_i(s),
\end{align*}
where
\[
C^n_i(t)=\sum_{k=1}^{A^n_0(t)}p_{n,\calR^n_i(s^n(k)-)}
=\int_{[0,t]}p_{n,\calR^n_i(s-)} dA^n_0(s),
\]
showing that $A^n_i-C^n_i=M^n_{i,1}$ is a martingale.

In the expression for $M^n_{i,2}$, the integrand is $\{\calF^n_t\}$-adapted
and has LCRL sample paths,
while the integrator is a martingale on the filtration. As a result,
$M^n_{i,2}$ is a local martingale
\cite[Theorem II.20]{protter}; using the estimate $\|M^n_{i,2}\|^*_t\le A^n_0(t)+c$
shows it is in fact a martingale. As a result, so is $\hat M^{A,n}_i$.
Finally, the expression for the quadratic variation is straightforward.
\qed

We will also need the following simple fact.

\begin{lemma}\label{lem5}
Let $M_N$, $N\in\Z_+$ be a martingale with $M_0=0$,
for which the increments $\Del_N=M_N-M_{N-1}$ satisfy
$\EE(|\Del_N|1_{\{|\Del_N|>a\}})\le\bar r(a)$, $a\ge0$, $N\in\N$,
and $\bar r(a)\to0$ as $a\to\iy$.
Then $N^{-1}\EE\|M\|^*_N<A(N)\to0$, where $\{A(N)\}$ depend only on
$\bar r$.
\end{lemma}

\proof
Let $b_N=\E[\Del_N1_{\{|\Del_N|\le a\}}]$ and note that
$b_N=-\E[\Del_N1_{\{|\Del_N|>a\}}]$ and consequently $|b_N|\le\bar r(a)$.
Write
\[
M_N=P_N+Q_N,
\qquad
P_N=\sum_{i=1}^N\{\Del_i1_{\{|\Del_i|>a\}}+b_i\},
\qquad
Q_N=\sum_{i=1}^N\{\Del_i1_{\{|\Del_i|\le a\}}-b_i\}.
\]
The quadratic variation
of the martingale $Q_N$ is bounded by $(a+\bar r(a))^2n$, giving
$\EE[\|Q\|^*_N]\le c(a+\bar r(a))n^{1/2}$, where $c^2$ is the constant from
the Burkholder-Davis-Gundy (BDG) inequality with $p=2$. For $P_N$,
\[
|P_N|\le \sum_{i=1}^N|\Del_i|1_{\{|\Del_i|>a\}}+N\bar r(a),
\]
thus $\E[\|P\|^*_N]\le2N\bar r(a)$. This gives
$n^{-1}\E[\|M\|^*_N]\le c(a+\bar r(a))N^{-1/2}+2\bar r(a)$.
Taking $a=N^{1/4}$ completes the proof.
\qed

\subsection{Uniform estimates on individual queue length processes}
\label{sec33}

The goal here is to calculate the dynamics of the individual
rescaled processes and develop estimates showing that they are
$C$-tight uniformly in $i\in[n]$.

In the following lemma, part (i) provides an equation, \eqref{10},
for the dynamics of individual queue lengths, and part (ii)
shows that, at the cost of introducing an error term,
one can replace the term $\hat S^n_i(T^n_i)$ in that equation
by a martingale. Both representations \eqref{10} and \eqref{22} are used
in this paper, where the former is used for estimates on each $\hat X^n_i$
that are uniform in $i$, and the latter is convenient for representations of
$\bar\xi^n$, as the martingale terms add up to a martingale.
The last part of the lemma uses \eqref{10} and gives uniform
second moment and tightness estimates.

\begin{lemma}\label{lem4}
i. One has
\begin{align}\label{10}
\hat X^n_i(t)&=\hat X^n_i(0-)+\hat E^n_i(t)+\hat A^n_i(t)
-\hat S^n_i(T^n_i(t))+\hat b^n_1t+\hat L^n_i(t).
\end{align}
Moreover, the sample paths of $\hat L^n_i$ are in $C^\up$ and
\begin{equation}\label{11}
\int_0^\iy\hat X^n_i(t)d\hat L^n_i(t)=0.
\end{equation}

\noi
ii. One has
\begin{equation}
\label{52}
\hat X^n_i(t)=\hat X^{1,n}_i(t)+e^{1,n}_i(t),
\end{equation}
where
\begin{equation}\label{22}
\hat X^{1,n}_i(t)=\hat X^n_i(0-)+\hat E^n_i(t)
+\hat A^n_i(t)-\hat M^{\dep,n}_i(t)
+\hat b^n_1t+\hat L^n_i(t),
\end{equation}
\[
e^{1,n}_i(t)=\tilde Z^n_i(0)-\tilde R^n_i(T^n_i(t))
=\tilde\mu^n(Z^n_i(0)-R^n_i(T^n_i(t))).
\]
iii. For $H^n_i=\hat S^n_i,\hat E^n_i,\hat A^n_i,\hat L^n_i$
and $\hat X^n_i$, one has
\begin{equation}\label{65}
\sup_n\sup_{i\in[n]}\EE[(\|H^n_i\|^*_t)^2]<\iy,\qquad t\ge 0,
\end{equation}
and for every $t$, $\eps>0$ and $\eta>0$ there is $\del>0$ such that
\begin{equation}\label{66}
\limsup_n\max_{i\in[n]}\PP(w_t(H^n_i,\del)>\eps)<\eta.
\end{equation}
\end{lemma}

\proof
i. By \eqref{01} and \eqref{07},
\begin{align*}
n^{-1/2}X^n_i(t)&=n^{-1/2}X^n_i(0-)+n^{-1/2}(E^n_i(t)-\la^nt)
+n^{-1/2}(\la^n-n\la)t+n^{1/2}\la t\\
&\quad+n^{-1/2}A^n_i(t)
-n^{-1/2}(S^n_i(T^n_i(t))-\mu^nT^n_i(t))-n^{-1/2}\mu^nT^n_i(t),
\end{align*}
and
\[
-n^{-1/2}\mu^nT^n_i(t)=-n^{-1/2}(\mu^n-n\mu)t-n^{1/2}\mu t+\hat L^n_i(t).
\]
Using \eqref{09}, \eqref{07+} and \eqref{08} gives \eqref{10}.
The properties of $\hat L^n_i$ and \eqref{11} follow from \eqref{03}.

ii.
Using \eqref{34} in \eqref{10} gives \eqref{52}.

iii.
By the central limit theorem for renewal processes \cite[\S 17]{bil},
and the fact that, by \eqref{04} and \eqref{06}, $n^{-1}(\la^n,\mu^n)\to(\la,\mu)$,
for each $i$, $(\hat E^n_i,\hat S^{0,n}_i)$
converge in law to $(E,S)$, a pair of mutually independent BM
starting at zero, with zero drift,
and diffusion coefficients $\la^{1/2}$ and $\mu^{1/2}\sig_\ser$, respectively
(where we recall $\la=\mu$).

We prove that \eqref{65} and \eqref{66} hold for
$\hat S^n_i$ by relating these processes to $\hat S^{0,n}_i$,
whose laws do not depend on $i$.
To this end, note that
\[
S^n_i(t)=S^{0,n}_i((t-Z^n_i(0))^+)-1_{\{t<Z^n_i(0)\}}
=S^{0,n}_i(t-Z^{\#,n}_i)-1_{\{t<Z^n_i(0)\}},
\]
where $Z^{\#,n}_i=Z^{\#,n}_i(t):=Z^n_i(0)\w t$.
Hence by \eqref{07+},
\[
\hat S^n_i(t)=\hat S^{0,n}_i(t-Z^{\#,n}_i)-\tilde\mu^n Z^{\#,n}_i
-n^{-1/2}1_{\{t<Z^n_i(0)\}}.
\]
As a result,
\begin{equation}\label{61}
\|\hat S^n_i\|^*_t \le \|\hat S^{0,n}_i\|^*_t+\tilde\mu^nZ^n_i(0)+n^{-1/2}
=\|\hat S^{0,n}_i\|^*_t+\tilde Z^n_i(0)+n^{-1/2},
\end{equation}
and
\[
w_t(\hat S^n_i,\del) \le w_t(\hat S^{0,n}_i,\del)+\tilde\mu^n Z^n_i(0)
=w_t(\hat S^{0,n}_i,\del)+\tilde Z^n_i(0).
\]
The latter inequality and the fact that the law of $\hat S^{0,n}_i$
does not depend on $i$ gives
\[
\max_{i\in[n]}\PP(w_t(\hat S^n_i,\del)>\eps)
\le \PP\Big(w_t(\hat S^{0,n}_1,\del)>\frac{\eps}{2}\Big)
+ \max_{i\in[n]}\PP\Big(\tilde Z^n_i(0)>\frac{\eps}{2}\Big).
\]
Using \eqref{40} and the $C$-tightness of $\hat S^{0,n}_1$, $n\in\N$,
shows that \eqref{66} is satisfied by $\hat S^n_i$.

Next, under our second moment assumptions on the inter-renewal times,
it is well known that the rescaled non-delayed renewal processes satisfy
\begin{equation}\label{303}
\sup_n\EE[(\|\hat S^{0,n}_1\|^*_t)^2]<\iy,
\end{equation}
for every $t$ \cite[Appendix 1]{krichagina1992diffusion}.
Since the law of $\hat S^{0,n}_i$ does not depend on $i$,
it follows from \eqref{61} and our assumption \eqref{62}
that $\hat S^n_i$ satisfy \eqref{65}.

It follows now that $\hat E^n_i$ satisfies both estimates,
for the law of $E^n_i$ is merely a special case of the law
of the non-delayed renewal processes $S^{0,n}_i$.

As for the processes $A^n_i$, recall from Lemma \ref{lem1}.iii that
$\hat M^{A,n}_i$ is a martingale and that
$[\hat M^{A,n}_i](t)=n^{-1}A^n_i(t)$.
One has
\[
\EE[n^{-1}A^n_i(t)]=n^{-1}\la^n_0\int_0^t\EE[p_{n,\calR^n_i(s)}]ds
\le cn^{-1}n^{3/2}n^{-1}t,
\]
where we used \eqref{17} and \eqref{05}.
Since this bound does not depend on $i$, it follows by the BDG
inequality that
\begin{equation}\label{18}
\lim_n\max_{i\in[n]}\EE[(\|\hat M^{A,n}_i\|^*_t)^2]=0.
\end{equation}
Again by the bound on $p_{n,r}$,
the Lipschitz constant of $\hat C^{A,n}_i$ is bounded
by $cn^{3/2}n^{-1/2}n^{-1}=c$.
It follows that both estimates \eqref{65} and \eqref{66} are satisfied
by $\hat A^n_i$.

To treat $\hat L^n_i$, recall Skorohod's lemma \cite[\S 8]{chu-wil},
stating that
for a trajectory $y\in D(\R_+,\R)$, if $(x,z)\in D(\R_+,\R_+)^2$ are such that
$x=y+z$, $z$ is nondecreasing and,
with the convention $z(0-)=0$, $\int_{[0,\iy)}xdz=0$,
then $z$ is given by
\[
z(t)=\sup_{s\in[0,t]}y^-(s),
\qquad
t\ge0.
\]
In particular,
\begin{equation}\label{12}
z(t)\le\|y\|_t,
\qquad
w_t(z,\del)\le w_t(y,\del),
\qquad
t>0,\,\del>0.
\end{equation}
In view of \eqref{10} and \eqref{11}, it follows that
\begin{equation}\label{16}
\hat L^n_i(t)\le\|\hat U^n_i\|_t,\qquad w_t(\hat L^n_i,\del)\le w_t(\hat U^n_i,\del),
\end{equation}
where
\begin{equation}\label{55}
\hat U^n_i(t)
=\hat X^n_i(0-)+\hat E^n_i(t)+\hat A^n_i(t)-\hat S^n_i(T^n_i(t))+\hat b^n_1t,
\end{equation}
and
\begin{equation}\label{67}
\hat X^n_i=\hat U^n_i+\hat L^n_i.
\end{equation}
Because $T^n_i$ are $1$-Lipschitz, we have that
$\hat S^n_i(T^n_i(\cdot))$ satisfy \eqref{65} and \eqref{66}.
Hence by the second moment assumption on initial condition
\eqref{53} and the boundedness of $\hat b^n_1$,
the same holds for $\hat U^n_i$. Finally, by \eqref{16}
and \eqref{67}, this is also true for $\hat L^n_i$ and $\hat X^n_i$.
\qed

\subsection{Empirical process tightness}\label{sec34}

\begin{lemma}\label{lem3}
$\bar\xi^n$ are $C$-tight in $D(\R_+,\calM_1)$.
\end{lemma}

\proof
Denote by $C^\eps$ the $\eps$-neighborhood, in $\R_+$, of a set $C\in\R_+$, and let
\[
d_{\rm L}(p,q)=\inf\{\eps>0: p(C^\eps)+\eps\ge q(C) \text{ and } q(C^\eps)+\eps\ge p(C)
\text{ for all } C\in\calB(\R_+)\},
\]
$p, q\in\calM_1$, denote the Levy-Prohorov metric, which
induces the topology of weak convergence on $\calM_1$.
Since $\bar\xi^n_0$ converge in probability, proving $C$-tightness of $\bar\xi^n$
amounts to showing that
for every $\eps>0$ and $\eta>0$ there exists $\del>0$ such that
\begin{equation}\label{53+}
\limsup_n\PP(w_T(\bar\xi^n,\del)>\eps)<\eta.
\end{equation}
Here and below we use the same notation $w_T$ for the metric $d_{\rm L}$ and for
the usual metric on $\R$.

We show that \eqref{53+} is a consequence of Lemma \ref{lem4}.iii.
Fix $\eps>0$ and $\eta>0$.
Given $n$ and $\del>0$, consider the event $\Om^n_\del=\{w_T(\bar\xi^n,\del)>\eps\}$.
On this event
there are $0\le s\le t\le T$, $t-s\le\del$, and $C\in\calB(\R_+)$, such that
\[
\xi^n_t(C^\eps)+n\eps < \xi^n_s(C)
\quad \text{ or } \quad
\xi^n_s(C^\eps)+n\eps < \xi^n_t(C).
\]
In both cases, the number of trajectories $\hat X^n_i$ whose displacement between
times $s$ and $t$ exceeds $\eps$ is greater than $n\eps$. Therefore
\[
\#\{i\in[n]:w_T(\hat X^n_i,\del)\ge\eps\}> n\eps.
\]
Hence by Chebychev's inequality,
\[
\PP(\Om^n_\del)\le(n\eps)^{-1}n\max_i\PP(w_T(\hat X^n_i,\del)>\eps).
\]
In view of \eqref{66}, given any $\eta_1>0$ there is $\del>0$ such that for all large $n$,
the maximum over $i$ in the above display is $<\eta_1$.
Hence, for such $\del$, the above is $\le \eps^{-1}\eta_1$.
Choosing $\eta_1$ such that $\eps^{-1}\eta_1<\eta$ and the corresponding $\del$
gives \eqref{53+}.
\qed

\section{Proof}
\label{sec4}

In \S\ref{sec41} we derive an equation for the empirical measure,
namely \eqref{24+} below, that is a precursor to the parabolic PDE.
That is, taking limits in \eqref{24+} gives rise to a weak solution to equation \eqref{70}.
Specifically, \S\ref{sec41} and, respectively, \S\ref{sec42},
show that the linear and nonlinear terms in \eqref{24+} converge to the corresponding
terms in \eqref{70}. In \S\ref{sec43}, \eqref{80} is obtained from \eqref{70}.
Based on this convergence, the proof of all the main results is then given
in \S\ref{sec44}.

\subsection{An equation for the empirical process}\label{sec41}

Fix a test function $\tilde\phi$ as in Definition~\ref{222},
that is, $\tilde\phi\in C^\iy_0(\R_+)$, $\tilde\phi(0)=0$.
Let
\begin{equation}\label{019}
\phi(x)=\int_0^x\tilde\phi(y)dy, \qquad x\in\R_+.
\end{equation}
Then $\phi\in C^\iy_b(\R_+)$ and $\phi(0)=\phi'(0)=0$.
Moreover, $\phi,\phi'$ and $\phi''$ are uniformly continuous on $\R_+$.
Apply $\phi$ to the dynamics \eqref{22}, noting that, unlike $X^n_i$,
$X^{1,n}_i$ may assume negative values.
On the r.h.s.\ of \eqref{22}, the terms $\hat A^n_i$ and $\hat M^{\dep,n}_i$
are piecewise constant. Thus
\begin{align}\label{23}
\notag
\phi(\hat X^{1,n}_i(t)) &= \phi(\hat X^n_i(0-))
+\int_0^t\phi'(\hat X^{1,n}_i(s))(d\hat E^{n,c}_i(s)+\hat b^n_1ds)
+e^{2,n}_i(t)
\\
&\quad
+\sum_{s\le t}(\phi(\hat X^{1,n}_i(s))-\phi(\hat X^{1,n}_i(s-))),
\end{align}
where $\hat E^{n,c}_i$ is the continuous part of $\hat E^n_i$,
\[
e^{2,n}_i(t)=\int_0^t(\phi'(\hat X^{1,n}_i(s))-\phi'(\hat X^n_i(s)))d\hat L^n_i(s),
\]
and we have used \eqref{11} and the continuity of the sample paths
of $\hat L^n_i$ to write
\[
\int\phi'(\hat X^n_i(s))d\hat L^n_i(s)=\phi'(0)\hat L^n_i(t)=0.
\]
In the last term of \eqref{23}, by Taylor's expansion,
jumps according to $\hat E^n_i$ can be expressed as
\[
\phi'(\hat X^{1,n}_i(s-))\Del\hat E^n_i(s)
+\frac{1}{2}\phi''(\hat X^{2,n}_i(s))\Del\hat E^n_i(s)^2,
\]
where
$\hat X^{2,n}_i(s)$ is an intermediate value between $\hat X^{1,n}_i(s-)$
and $\hat X^{1,n}_i(s)$
(we leave unspecified the value of the processes $\hat X^{2,n}_i$
away from times of jumps of $\hat X^{1,n}_i$).
Similarly, jumps according to $\hat M^{\dep,n}_i$ are expressed as
\[
-\phi'(\hat X^{1,n}_i(s-))\Del\hat M^{\dep,n}_i(s)
+\frac{1}{2}\phi''(\hat X^{2,n}_i(s))\Del\hat M^{\dep,n}_i(s)^2,
\]
where again $\hat X^{2,n}_i(s)$ are intermediate points,
and because, a.s., jumps of $\hat E^n_i$ and $\hat M^{\dep,n}_i$
do not occur simultaneously (for the same $i$),
we may express the intermediate values by the same process.
Finally, jumps according to $\hat A^n_i$ -- here we only need first order approximation -- are expressed as
\[
\phi'(\hat X^{2,n}_i(s))\Del\hat A^n_i(s).
\]
Note that jumps of $\hat E^n_i$ and $\hat A^n_0$ are of size
$n^{-1/2}$. Hence in all three cases,
\begin{equation}\label{24}
|\hat X^{2,n}_i(t)-\hat X^{1,n}_i(t-)|\le |\Del\hat X^{1,n}_i(t)|
\le \max\{|\Del\hat M^{\dep,n}_i(t)|, n^{-1/2}\},
\qquad t\in\calJ^n_i,
\end{equation}
where $\calJ^n_i$ is the set of jump times of $\hat X^{1,n}_i$.

By \eqref{23}, using $\hat E^n_i=\hat E^{n,c}_i+\Del\hat E^n_i$ and
\[
\int_0^t\phi'(\hat X^{1,n}_i(s))d\hat E^{n,c}_i(s)
=\int_0^t\phi'(\hat X^{1,n}_i(s-))d\hat E^{n,c}_i(s),
\]
we have
\begin{align}\label{24+}
\lan \phi,\bar\xi^n_t\ran
=\lan \phi,\bar\xi^n_{0-}\ran
+\int_0^t\lan b_1\phi'+\frac{\sig^2}{2}\phi'',\bar\xi^n_s\ran ds
+\Gam^n(t)+\sum_{j=1}^6f^{j,n}(t),
\end{align}
where
\begin{equation}\label{57}
\Gam^n(t)=\frac{b_0}{n}\sum_i\int_0^t\phi'(\hat X^n_i(s))
\Big(\frac{n-\calR^n_i(s)}{n}\Big)^{\ell-1}ds
\end{equation}
is the interaction term,
\begin{align*}
f^{1,n}(t)&=f^{1,n}_1(t)+f^{1,n}_2(t)\\
&:=
\frac{1}{n}\sum_i\int_{[0,t]}\phi'(\hat X^{1,n}_i(s-))d\hat E^n_i(s)
-\frac{1}{n}\sum_i\int_{[0,t]}\phi'(\hat X^{1,n}_i(s-))d\hat M^{\dep,n}_i(s),
\\
f^{2,n}(t)&=\lan \phi,\bar\xi^n_t\ran-\lan \phi,\bar\xi^{1,n}_t\ran,
\\
f^{3,n}(t)&=\frac{1}{n}\sum_ie^{2,n}_i(t),
\\
f^{4,n}(t)&=\hat b^n_1\int_0^t\lan \phi',\bar\xi^{1,n}_s\ran ds
-b_1\int_0^t\lan \phi',\bar\xi^n_s\ran ds,
\\
f^{5,n}(t)
&=\frac{1}{2n}\sum_i\int_{[0,t]}\phi''(\hat X^{2,n}_i(s))
\{d[\hat E^n_i](s)+d[\hat M^{\dep,n}_i](s)\}
-\frac{\sig^2}{2}\int_0^t\lan \phi'',\bar\xi^n_s\ran ds,
\\
f^{6,n}(t)
&=
\frac{1}{n}\sum_i\int_{[0,t]}\phi'(\hat X^{2,n}_i(s))d\hat A^n_i(s)
-\Gam^n(t)
\end{align*}
are ``error'' terms, and
\[
\bar\xi^{1,n}_t=n^{-1}\sum_{i\in[n]}\del_{\hat X^{n,1}_i(t)}.
\]
Note that the terms $f^{5,n}$ and $f^{6,n}$, which involve $\hat X^{2,n}_i$,
are well defined, for their evaluation requires
the values of $\hat X^{2,n}_i(t)$ only at $\calJ^n_i$.

\begin{lemma}\label{lem2}
With $\phi$ as in \eqref{019}, for $1\le j\le 6$, $f^{j,n}\to0$ in probability
in $D(\R_+,\R)$.
\end{lemma}

We note that one of the reasons the proof is somewhat involved
is that we work under minimal moment assumptions.
For example, the second order term $f^{5,n}$ involves
martingales given in terms of squares of the primitive data $\zeta^n_i(k)$.
These martingales do not, in general, possess moments higher than first
because only the second moments of $\zeta^n_i(k)$ are assumed finite.
Hence estimates on these martingales based on
their quadratic variation are not available. Our treatment of this
difficulty is based on Lemma \ref{lem5}.

\proof
\noi{\it Step 1.}
Estimating $f^{1,n}$.
The cross variation between independent Poisson processes is zero,
hence $[\hat E^n_i,\hat E^n_j](t)=0$ for $i\ne j$.
Hence by Lemma \ref{lem1},
\[
[f^{1,n}_1](t)\le cn^{-2}\sum_i[\hat E^n_i](t)
\le cn^{-3}\sum_iE^n_i(t).
\]
Since each $E^n_i$ is a Poisson process of intensity $\le cn$,
$\EE\{[f^{1,n}_1](t)\}\le ctn^{-1}$, and it follows by the BDG inequality
that $f^{1,n}_1\to0$ in probability.
As for $f^{1,n}_2$, we have by Lemma \ref{lem1}.ii that
$[\hat M^{\dep,n}_i,\hat M^{\dep,n}_j](t)$,
$i\ne j$, has zero mean. Hence by Lemma \ref{lem1}.i,
\[
\EE\{[f^{1,n}_2](t)\}\le cn^{-2}\sum_i\EE\{[\hat M^{\dep,n}_i](t)\}
\le cn^{-3}\sum_i\EE\{D^n_i(t)\}\le cn^{-1}t.
\]
Thus $f^{1,n}_2\to0$ in probability,
and we conclude that $f^{1,n}\to0$ in probability.

\noi{\it Step 2.}
We show that for all $t>0$ and $\eps>0$,
\begin{equation}\label{50}
\lim_{n\to\iy}\max_{i\in[n]}\PP(\|e^{1,n}_i\|^*_t>\eps)=0.
\end{equation}
In accordance with \eqref{40} we denote
$\tilde Z^n_i(k)=\tilde\mu^nZ^n_i(k)$, $k\ge1$.
By the definition of $e^{1,n}_i(t)$ and $T^n_i(t)\le t$,
\begin{align}
\notag
|e^{1,n}_i(t)|\le\tilde Z^n_i(0)\vee\|\tilde R^n_i\|^*_t
&\le \max\{\tilde Z^n_i(k): 0\le k\le S^n_i(t)\}
\\
\label{42}
&\le \max\{\tilde Z^n_i(k): 0\le k\le S^{0,n}_i(t)\}.
\end{align}
Thus, by \eqref{40} and the fact that the law of
$((\tilde Z^n_i(k),k\ge1),S^{0,n}_i(t))$ does not depend on $i$,
it suffices to prove that
\[
\tilde Y^n:=\max\{\tilde Z^n_1(k): 1\le k\le S^{0,n}_1(t)\}\to0
\text{ in probability.}
\]
This can be argued via the $C$-tightness of $\{\hat S^{0,n}_1\}$ as follows.
Given $n$ and $t_1<t_2$, we have
\begin{equation}\label{41}
\text{if $S^{0,n}_1(t_1)=S^{0,n}_1(t_2)$ then }
\hat S^{0,n}_1(t_1)-\hat S^{0,n}_1(t_2)=\tilde\mu^n(t_2-t_1).
\end{equation}
This shows that
\[
\max\{Z^n_1(k): 1\le k\le S^{0,n}_1(t)-1\}
\le 2(\tilde\mu^n)^{-1}\|\hat S^{0,n}_1\|^*_t.
\]
To include also $k=S^{0,n}_1(t)$, let the residual time process
$R^{0,n}_1$ be defined
analogously to $R^n_1$, for $S^{0,n}_1$ in place of $S^n_1$.
Note that if $R^{0,n}_1(t)>1$ holds then \eqref{41} holds
with $t_1=t$ and $t_2=t+1$ hence
$\hat S^{0,n}_1(t)-\hat S^{0,n}_1(t+1)=\tilde\mu^n$.
Because $\tilde\mu^n\to\iy$
and $\|\hat S^{0,n}_1\|^*_{t+1}$ are tight, this shows that
w.h.p., $R^{0,n}_1(t)\le 1$. As a result, again by \eqref{41}, w.h.p.,
\[
Y^n:=\max\{Z^n_1(k): 1\le k\le S^{0,n}_1(t)\}
\le 2(\tilde\mu^n)^{-1}\|\hat S^{0,n}_1\|^*_{t+1}.
\]
This shows that $Y^n\to0$ in probability.
Next, given $\del>0$, on the event
$Y^n<\del$, we have by \eqref{41},
\begin{align*}
w_{t+1}(\hat S^{0,n}_1,\del)
&\ge
\sup\{|\hat S^{0,n}_1(t_1)-\hat S^{0,n}_1(t_2)|:0\le t_1<t_2\le t+1,
\, S^{0,n}_1(t_1)=S^{0,n}_1(t_2)\}
\\
&\ge
\tilde\mu^n\sup\{t_2-t_1:0\le t_1<t_2\le t+1,
\, S^{0,n}_1(t_1)=S^{0,n}_1(t_2)\}
\\
&\ge\tilde\mu^nY^n.
\end{align*}
Hence, w.h.p., $\tilde Y^n=\tilde\mu^nY^n\le w_{t+1}(\hat S^{0,n}_1,\del)$.
Using the $C$-tightness of $\hat S^{0,n}_1$, sending
$n\to\iy$ and then $\del\to0$ shows that $\tilde Y^n\to0$ in probability.
This proves \eqref{50}.

\noi{\it Step 3.}
We can now control $f^{2,n}$ and $f^{4,n}$. Recall \eqref{52}. Then,
denoting by $m_\phi(\cdot)$ the modulus of continuity of $\phi$, for any $\del>0$,
\[
\|f^{2,n}\|^*_t\le\frac{1}{n}\sum_i\|\phi(\hat X^n_i)-\phi(\hat X^{1,n}_i)\|^*_t
\le m_\phi(\del)+\frac{2\|\phi\|_\iy}{n}\#\{i:\|e^{1,n}_i\|^*_t>\del\}.
\]
Given $\eps>0$ let $\del>0$ be such that $m_\phi(\del)<\eps/2$. Then
\[
\PP(\|f^{2,n}\|^*_t>\eps)\le4\|\phi\|_\iy\eps^{-1}\max_i\PP(\|e^{1,n}_i\|^*_t>\del).
\]
Sending $n\to\iy$, the above converges to $0$ by \eqref{50},
proving that $f^{2,n}\to0$ in probability.

Next, similarly, $\lan \phi',\bar\xi^{1,n}\ran-\lan \phi',\bar\xi^n\ran\to0$
in probability, and
\[
\|f^{4,n}\|^*_t\le
t\,|\hat b^n_1|\, \|\lan \phi',\bar\xi^{1,n}\ran-\lan \phi',\bar\xi^n\ran\|^*_t
+
|\hat b^n_1-b_1|\int_0^\cdot\lan \phi',\bar\xi^n_s\ran ds.
\]
The first term on the right converges to $0$ in probability
by the argument shown for $f^{2,n}$ and the boundedness of $\hat b^n_1$.
In the second term, the integral is bounded by $t\|\phi'\|_\iy$,
hence the convergence of this term to zero
follows from $\hat b^n_1\to b_1$. This proves $f^{4,n}\to0$ in probability.

\noi{\it Step 4.} To estimate $f^{3,n}$, note that, for any $\del>0$,
\[
\|f^{3,n}\|^*_t\le\frac{m_{\phi'}(\del)}{n}\sum_i\hat L^n_i(t)
+\frac{2\|\phi'\|_\iy}{n}\sum_i 1_{\{\|e^{1,n}_i\|^*_t>\del\}}\hat L^n_i(t).
\]
Hence
\[
\EE[\|f^{3,n}\|^*_t]\le\frac{m_{\phi'}(\del)}{n}\sum_i\EE[\hat L^n_i(t)]
+\frac{2\|\phi'\|_\iy}{n}\sum_i \PP(\|e^{1,n}_i\|^*_t>\del)^{1/2}
(\EE[\hat L^n_i(t)^2])^{1/2}.
\]
Since by Lemma \ref{lem4}.iii, $\max_{i\in[n]}\EE[\hat L^n_i(t)^2]<c$, it follows that
\[
\EE[\|f^{3,n}\|^*_t]\le cm_{\phi'}(\del)+c\max_{i\in[n]}\PP(\|e^{1,n}_i\|^*_t>\del).
\]
In view of \eqref{50}, if we take $n\to\iy$ and then $\del\to0$,
the expression on the right converges to $0$,
hence $f^{3,n}\to0$ in probability.

\noi{\it Step 5.}
Recall
\[
f^{5,n}(t)
=\frac{1}{2n}\sum_i\int_{[0,t]}\phi''(\hat X^{2,n}_i(s))
\{d[\hat E^n_i](s)+d[\hat M^{\dep,n}_i](s)\}
-\frac{\sig^2}{2}\int_0^t\lan\phi'',\bar\xi^n_s\ran ds,
\]
\[
[\hat E^n_i](t)=n^{-1}E^n_i(t)
\qquad
[\hat M^{\dep,n}_i](t)=n^{-1}\sum_{k=1}^{D^n_i(t)}\zeta^n_i(k)^2.
\]
Let
\begin{align*}
f^{5,n}_1&=\frac{1}{n}\sum_i\int_{[0,t]}\phi''(\hat X^{1,n}_i(s-))d[\hat E^n_i](s)
-\la\int_0^t\lan\phi'',\bar\xi^{1,n}_s\ran ds,
\\
f^{5,n}_2&=\frac{1}{n}\sum_i\int_{[0,t]}\phi''(\hat X^{1,n}_i(s-))d[\hat M^{\dep,n}_i](s)
-\la\sig_\ser^2\int_0^t\lan\phi'',\bar\xi^{1,n}_s\ran ds,
\\
f^{5,n}_3&=\int_0^t\lan\phi'',\bar\xi^n_s\ran ds-\int_0^t\lan\phi'',\bar\xi^{1,n}_s\ran ds
\\
f^{5,n}_4&=
\frac{1}{n}\sum_i\int_{[0,t]}(\phi''(\hat X^{2,n}_i(s))-\phi''(\hat X^{1,n}_i(s-)))
\{d[\hat E^n_i](s)+d[\hat M^{\dep,n}_i](s)\}.
\end{align*}
Then $f^{5,n}=\sum_{j}f^{5,n}_j$.
For $f^{5,n}_1$, write
\begin{align*}
f^{5,n}_1&=\frac{1}{n}\sum_i\int_{[0,t]}\phi''(\hat X^{1,n}_i(s-))n^{-1}dE^n_i(s)
-\la\int_0^t\lan\phi'',\bar\xi^{1,n}_s\ran ds
\\
&=\frac{1}{n}\sum_i\int_{[0,t]}\phi''(\hat X^{1,n}_i(s-))n^{-1/2}d\hat E^n_i(s)
+\frac{1}{n}\sum_i\int_0^t\phi''(\hat X^{1,n}_i(s))n^{-1/2}\hat\la^nds.
\end{align*}
In the first sum, each term is a martingale having expected quadratic variation bounded by
$c n^{-2}\la^nt\le cn^{-1}t$
where we use Lemma \ref{lem1}.i, and $c$ does not depend on $i$ or $n$.
Hence the first normalized sum is a martingale whose expected quadratic variation
is $\le cn^{-2}t$. In particular, it converges to zero in probability.
In the second sum, each term is bounded in absolute value by $ctn^{-1/2}$,
hence the second normalized sum converges to zero. This shows
$f^{5,n}_1\to0$ in probability.

Using Lemma \ref{lem1} for the expression for $[\hat M^{\dep,n}_i]$
and denoting
\[
P^n_i(t)=\phi''(\hat X^{1,n}_i(t)),
\qquad
Q^n_i(t)=\sum_{k=1}^{D^n_i(t)}q^n_i(k),
\qquad
q^n_i(k)=\zeta^n_i(k)^2-\sig_\ser^2,
\]
write $f^{5,n}_2$ as
\begin{align}\label{091}
f^{5,n}_2&=\frac{1}{n}\sum_i\int_{[0,t]}P^n_i(s-)n^{-1}dQ^n_i(s)
+\frac{\sig_\ser^2}{n}\sum_i\int_{[0,t]}P^n_i(s-)\{n^{-1}dD^n_i(s)
-\la ds\}.
\end{align}
To estimate the first term above, note that, by \eqref{07+},
\[
D^n_i(t)= S^n_i(T^n_i(t))\le S^n_i(t)=n^{1/2}\hat S^n_i(t)+\mu^nt
\le n^{1/2}\hat S^n_i(t)+c_1nt,
\]
for a constant $c_1$. For any $\beta>2c_1t$,
\[
\max_{i\in[n]}\PP(D^n_i(t)>\beta n)
\le\max_{i\in[n]}\PP(\hat S^n_i(t)>n^{1/2}\beta/2)\le cn^{-1}\beta^{-2},
\]
by Lemma \ref{lem4}.iii, where $c=c(t)$.
Denote by $0\le t^n_i(1)<t^n_i(2)<\cdots$ the jump times of $D^n_i$
and let $m^n_i(N)=\sum_{k=1}^NP^n_i(t^n_i(k)-)q^n_i(k)$.
Then the first term in \eqref{091} is given by
$n^{-2}\sum_im^n_i(D^n_i(t))$, which we write as
\begin{equation}\label{092}
\frac{1}{n^2}\sum_i m^n_i(D^n_i(t)\w(\beta n))
+
\frac{1}{n^2}\sum_i\{m^n_i(D^n_i(t)-m^n_i(D^n_i(t)\w(\beta n))\}.
\end{equation}
Now, $m^n_i(\cdot)$ is a discrete parameter martingale with
increments bounded by $\|\phi''\|_\iy |q^n_i(k)|$. Moreover, $q^n_i(k)$ are
IID in $k$ and $i$,
and possesses a first moment in view of our second moment
assumption on the service times. Hence Lemma \ref{lem5} is applicable,
showing that
\[
N^{-1}\EE\|m^n_i\|^*_N<A(N)\to0 \qquad \text{as } N\to\iy,
\]
where $A(\cdot)$ does not depend on $n$ ot $i$. Hence
\[
\EE\Big|\frac{1}{n^2}\sum_i m^n_i(D^n_i(t)\w(\beta n))\Big|
\le\frac{1}{n}\sum_in^{-1}\EE[\|m^n_i\|^*_{\beta n}]
\le n^{-1}\beta nA(\beta n) = \beta A(\beta n),
\]
where throughout this paragraph
$\beta n$ should be read as $\lfloor \beta n\rfloor$.
Hence the limit of the above expression is $0$ for any $\beta$.
Next, using $D^n_i\le S^n_i$,
\begin{align*}
\E\Big|\frac{1}{n^2}\sum_i\{m^n_i(D^n_i(t)-m^n_i(D^n_i(t)\w(\beta n))\}\Big|
&\le \|\phi''\|_\iy n^{-2}\sum_i 1_{\{S^n_i(t)>\beta n\}}\sum_{k=\beta n+1}^{S^n_i(t)}|q^n_i(k)|
\\
&=cn^{-2}\sum_i\sum_{k=1}^{(S^n_i(t)-\beta n)^+}|q^n_i(k-\beta n)|
\\
&=
cn^{-2}\sum_i\E[(S^n_i(t)-\beta n)^+]\E[|q^n_i(1)|]
\\
&\le cn^{-1}\max_{i\in[n]}\E[(S^n_i(t)-\beta n)^+],
\end{align*}
where Wald's identity is used on the third line.
Since it follows from Lemma \ref{lem4}.iii that $n^{-1}S^n_i(t)$
are uniformly integrable
in $i$ and $n$,
\[
\lim_{\beta\to\iy}\limsup_nn^{-1}\max_{i\in[n]}\E[(S^n_i(t)-\beta n)^+]=0.
\]
As a result, the expression in \eqref{092} converges
to zero in probability.

As for the second term in \eqref{091}. In view of Lemma \ref{lem4}.iii and \eqref{04}, one has
\[
\max_{i\in[n]}\|n^{-1}S^n_i-\la\iota\|^*_t\to0 \qquad \text{in probability}.
\]
Also from Lemma \ref{lem4}.iii and the relation \eqref{07} between
$\hat L^n_i$ and $T^n_i$, $\max_{i\in[n]}\|T^n_i-\iota\|^*_t\to0$ in probability.
Hence $D^n_i=S^n_i(T^n_i)$ also satisfies
$\kap_n:=\max_{i\in[n]}\|n^{-1}D^n_i-\la\iota\|^*_t\to0$ in probability.
Denote $\eps=t/j$, $s_k=k\eps$ and $\tilde P^n_i(s)=P^n_i(s-)$.
Then
\begin{align*}
\int_{(0,t]}\tilde P^n_i(s)\{n^{-1}dD^n_i(s)-\la ds\}
&=\sum_{k=0}^{j-1}\Big\{\int_{(s_k,s_{k+1}]}(\tilde P^n_i(s)-\tilde P^n_i(s_k))\{n^{-1}dD^n_i(s)-\la ds\}
\\
&\qquad\qquad
+\tilde P^n_i(s_k)\{n^{-1}D^n_i(s_{k+1})-n^{-1}D^n_i(s_k)-\la\eps\}\Big\}.
\end{align*}
This gives
\begin{align}\label{090}
\Big|\int_{(0,t]}\tilde P^n_i(s)\{n^{-1}dD^n_i(s)-\la ds\}\Big|
&\le w_T(\tilde P^n_i,\eps)(n^{-1}D^n_i(t)+\la t)
+\|\tilde P^n_i\|^*_t2j\kappa_n
\\
\notag
&\le w_T(\tilde P^n_i,\eps)(2\la t+\kappa_n)+\|\phi''\|_\iy 2j\kappa_n.
\end{align}
We recall from Lemma \ref{lem4}.iii that $\hat X^n_i$ are uniformly $C$-tight.
Along with the uniform estimate \eqref{50} on $e^{1,n}_i$, this shows
that so are $\tilde P^n_i$. In particular, these processes
satisfy \eqref{66}. Hence, on sending $n\to\iy$ and then $\eps\to0$,
it follows that the right side of \eqref{090} converges to zero in probability,
uniformly in $i\in[n]$. Hence the second term in $f^{5,n}_2$ converges
to zero in probability, and we conclude that $f^{5,n}_2\to0$ in probability.

Next, recalling the result regarding $f^{2,n}$ and replacing $\phi$ by $\phi''$ shows that
$f^{5,n}_3\to0$ in probability.

To bound $f^{5,n}_4$, fix $\eps>0$ and let $\del>0$ be such that $|x-y|<\del$ implies
$|\phi''(x)-\phi''(y)|<\eps$. Let
$\calJ^{n,\del}_i=\{s\in\calJ^n_i:|\hat X^{2,n}_i(s)-\hat X^{1,n}_i(s-)|>\del\}$. For $n$ sufficiently large,
the size of jumps of $\hat E^n$ (which is $n^{-1/2}$) is smaller than $\del$,
hence, by \eqref{24}, the corresponding jump times are not members of $\calJ^{n,\del}_i$. Thus
the $i$-th term in $f^{5,n}_4$ is bounded, in absolute value, by
\begin{align*}
&\eps([\hat E^n_i](t)+[\hat M^{\dep,n}_i](t))
+2\|\phi''\|_\iy\sum_{s\in\calJ^{n,\del}_i,s\le t}\Del[\hat M^{\dep,n}_i](s)
\\
&\le\eps n^{-1}\Big(E^n_i(t)+\sum_{k=1}^{D^n_i(t)}\zeta^n_i(k)^2\Big)
+cn^{-1}\sum_{k=1}^{D^n_i(t)}\zeta^n_i(k)^21_{\{n^{-1/2}|\zeta^n_i(k)|>\del\}}
\end{align*}
where we used the fact that for $s\in \calJ^{n,\del}_i$, $|\Del\hat X^{1,n}_i(s)|=|\Del\hat M^{\dep,n}_i(s)|$
by \eqref{22}. Using \eqref{64}, the expected value of the above expression is
\[
\le \eps n^{-1}(\la^n_i t+\sig_\ser^2\EE[D^n_i(t)])
+cn^{-1}\EE[\zeta^n_1(1)^21_{\{n^{-1/2}|\zeta^n_1(1)|>\del\}}]\EE[D^n_i(t)],
\]
where for the last term we used Wald's identity in exactly the same way
as in the proof of Lemma \ref{lem1}. In view of Lemma \ref{lem4}.iii,
$\E[S^n_i(t)]<cn$, where $c$ does not depend on $i$ or $n$. Since $D^n_i(t)\le S^n_i(t)$,
it follows that
\[
\E[|f^{5,n}_4|]\le c\eps+c\EE[\zeta^n_1(1)^21_{\{n^{-1/2}|\zeta^n_1(1)|>\del\}}].
\]
Taking $n\to\iy$ then $\eps\to0$ shows that $f^{5,n}_4\to0$ in probability.

\noi{\it Step 6.}
Finally,
\[
f^{6,n}(t)
=\frac{1}{n}\sum_i\int_{[0,t]}[\phi'(\hat X^{2,n}_i(s))-\phi'(\hat X^n_i(s))]d\hat A^n_i(s)
+\frac{1}{n}\sum_i\int_{[0,t]}\phi'(\hat X^n_i(s))d\hat A^n_i(s)
-\Gam^n(t).
\]
By Lemma \ref{lem1}.iii,
\begin{align*}
\frac{1}{n}\sum_i\int_{[0,t]}\phi'(\hat X^n_i(s-))d\hat A^n_i(s)
&=\la^n_0n^{-3/2}\sum_i\int_0^t\phi'(\hat X^n_i(s))p_{n,\calR^n_i(s)}ds
+\bar M^{A,n}(t),
\\
\bar M^{A,n}(t)
&=
\frac{1}{n}\sum_i\int_{[0,t]}\phi'(\hat X^n_i(s-))d\hat M^{A,n}_i(s).
\end{align*}
By Lemma \ref{lem1}.iii, $[\bar M^{A,n}](t)\le cn^{-3}\sum_iA^n_i(t)=cn^{-3}A^n_0(t)$.
Since $A^n_0$ is a Poisson process of intensity $\la^n_0\le cn^{3/2}$ (by \eqref{05}),
it follows that $\bar M^{A,n}\to0$ in probability.

Now, writing
\[
p_{n,r}=\frac{\ell}{n}\,\frac{(n-r)(n-r-1)\cdots(n-r-\ell+2)}{(n-1)(n-2)\cdots(n-\ell+1)}
\]
shows that
\begin{equation}\label{097}
p_{n,r}=\frac{\ell}{n}\Big[\Big(1-\frac{r}{n}\Big)^{\ell-1}+\al_{n,r}\Big],
\qquad
\al_n:=\max_{r\in[n]}|\al_{n,r}|\to0.
\end{equation}
Hence
\begin{align*}
&\la^n_0n^{-3/2}\sum_i\int_0^t\phi'(\hat X^n_i(s))p_{n,\calR^n_i(s)}ds-\Gam^n(t)
\\
&\qquad=
\sum_i\int_0^t\phi'(\hat X^n_i(s))\Big[\la^n_0n^{-3/2}p_{n,\calR^n_i(s)}
-\frac{b_0}{n}\Big(\frac{n-\calR^n_i(s)}{n}\Big)^{\ell-1}\Big]ds
\\
&\qquad=
\sum_i\int_0^t\phi'(\hat X^n_i(s))
\Big(\frac{n-\calR^n_i(s)}{n}\Big)^{\ell-1}
\Big[\la^n_0n^{-3/2}\frac{\ell}{n}-\frac{b_0}{n}\Big]ds
\\
&\qquad\qquad+\sum_i\int_0^t\phi'(\hat X^n_i(s))\la^n_0n^{-3/2}\frac{\ell}{n}\al_{n,\calR^n_i(s)}ds.
\end{align*}
We have $\la^n_0n^{-3/2}\ell\to b\ell=b_0$, which shows that the first
sum converges to $0$. From the fact $\al_n\to0$ we also have that
the last sum converges to $0$.

To bound the first term in $f^{6,n}$, note that $\hat A^n_i$ are nondecreasing.
Also, since the jumps of $\hat A^n_i$ are of size $n^{-1/2}$,
one has $|\phi'(\hat X^{2,n}_i(s)-\phi'(\hat X^{1,n}_i(s))|\le m_{\phi'}(n^{-1/2})$
at any jump time $s$ of $\hat A^n_i$. Moreover,
$\max_{i\in[n]}|\phi'(\hat X^{1,n}_i(s)-\phi'(\hat X^n_i(s))|
\le \max_{i\in[n]}m_{\phi'}(\|e^{1,n}_i\|^*_t)\to0$ in probability, by \eqref{50}.
Hence
\[
\max_{i\in[n]}\Big|\int_{[0,t]}[\phi'(\hat X^{2,n}_i(s))-\phi'(\hat X^n_i(s))]d\hat A^n_i(s)\Big|
\le \max_{i\in[n]}\{(m_{\phi'}(n^{-1/2})+m_{\phi'}(\|e^{1,n}_i\|^*_t))\hat A^n_i(t)\}
\]
converges to $0$
in probability, where Lemma \ref{lem4}.iii is used for a uniform estimate
on $\hat A^n_i(t)$.
It follows that $f^{6,n}\to0$ in probability.
\qed

\subsection{Interaction term under limit}\label{sec42}

In this subsection we prove the following.
\begin{lemma}\label{lem6}
Given a subsequence along which $\bar\xi^n\To\xi$, one has
$\Gam^n\To\Gam$ in $C(\R_+,\R)$ along the subsequence, where
\[
\Gam(t)=\frac{b_0}{\ell}\int_0^t\int_{\R_+}\phi'(x)\frS(\xi_s[x,\iy),\xi_s(x,\iy))\xi_s(dx)ds.
\]
\end{lemma}

Recall that if $\xi_t$ is atomless for all $t>0$ then the integrand
in the above expression simplifies to $\ell \phi'(x)\xi_s[x,\iy)^{\ell-1}$,
in which case $\Gam$ is directly related to the form of the interaction
term in equations \eqref{14_}, \eqref{100} and \eqref{111}.
However, the atomless property has not been established at this stage of the proof.

\proof
Invoking the Skorohod representation theorem, assume without loss
of generality that $\bar\xi^n\to\xi$ a.s.
Because $0\le\calR^n_i\le n$, the integrands in \eqref{57}
are bounded by $\|\phi'\|_\iy$. Hence by bounded convergence,
it suffices to prove that, a.s., for every $t$,
\begin{equation}\label{31}
\gamma^n(t):=
\frac{1}{n}\sum_i\phi'(\hat X^n_i(t))\bar\calR^{c,n}_i(t)^{\ell-1}
\to\frac{1}{\ell}\int_{\R_+}\phi'(x)\frS(\xi_t[x,\iy),\xi_t(x,\iy))\xi_t(dx),
\end{equation}
where
\[
\calR^{c,n}_i(t)=n-\calR^n_i(t),
\qquad
\bar\calR^{c,n}_i(t)=n^{-1}\calR^{c,n}_i(t).
\]
Fix $t$ and $\eps>0$.
The function $x\mapsto\xi_t(x,\iy)$
has at most countably many discontinuities. Hence one can find
a finite sequence $x_*=y_0<y_1<\ldots<y_K=x^*$ such that
$[x_*,x^*]$ contains the compact support of $\phi'$ (recall \eqref{019}),
$y_k-y_{k-1}<\eps$, $k=1,\ldots,K$, and moreover $\xi_t(\{y_k\})=0$, $k=0,\ldots,K$.
Because $y_k$ are not charged,
both $\bar\xi^n_t(y_{k-1},y_k)$ and $\bar\xi^n_t[y_{k-1},y_k]$
converge a.s.\ to $\xi_t(y_{k-1},y_k)$.
Define by $\rank^c(i;x)=n-\rank(i;x)$, $x\in\R^n$, $i\in[n]$
the complementary rank. Then $0\le\rank^c(i;x)\le n-1$ and by \eqref{17-},
\[
\rank^c(i;x)=\#\{j:x_j>x_i\}+\#\{j>i:x_j=x_i\}.
\]
Denote $x_i=\hat X^n_i(t)$, $x=(x_1,\ldots,x_n)$, and for $k=1,\ldots,K$,
\[
I_k=[y_{k-1},y_k),
\qquad
V_k=\{i:x_i\in I_k\},
\qquad
P_k=\xi^n_t(I_k),\qquad Q_k=\xi^n_t[y_k,\iy).
\]
Fix $k$, and for $i\in V_k$ let $j(i)=\rank(i;(x_l)_{l\in V_k})$.
This is a relabeling of $V_k$ according to the rank of its members within $V_k$.
With this notation we have
\[
\rank^c(i;x)=Q_k+P_k-j(i), \qquad i\in V_k.
\]
Because $j(i)$ take all values between $1$ and $P_k$ as $i$ varies
in $V_k$, we have
\[
\sum_{i\in V_k}\rank^c(i;x)^{\ell-1}
=\sum_{j=1}^{P_k}(Q_k+P_k-j)^{\ell-1}
=\sum_{j=0}^{P_k-1}(Q_k+j)^{\ell-1}.
\]
If $p_n$ is a sequence s.t.\ $n^{-1}p_n\to p\ge0$
then $n^{-1}\sum_{j=0}^{p_n}(j/n)^{\ell-1}\to\int_0^pz^{\ell-1}dz
=\ell^{-1}p^\ell$. Noting, as mentioned, that $n^{-1}P_k\to\xi_t(I_k)$,
we get
\begin{equation}
\label{60}
\frac{1}{n}\sum_{i\in V_k}\bar\calR^{c,n}_i(t)^{\ell-1}
\to \frac{1}{\ell}(\xi_t(y_{k-1},\iy)^\ell-\xi_t(y_k,\iy)^\ell)
=\frac{1}{\ell}\xi_t(I_k)\frS(\xi_t(y_{k-1},\iy),\xi_t(y_k,\iy)).
\end{equation}
Fix a sequence $\eps_m>0$, $\eps_m\to0$. For each $m$, let points as the above
$y_k$, corresponding to $\eps_m$, be denoted by $y_k^m$.
For $x\in[x_*,x^*]$ let $x_{(m)}$ and $x^{(m)}$ be the unique two points
$y_{k-1}^m$ and $y_k^m$, respectively, for which $x\in[y_{k-1}^m,y_k^m)$.
This gives $x^{(m)}\to x$, $x^{(m)}>x$. By right-continuity of
$x\mapsto\xi_t(x,\iy)$, it follows that $\xi_t(x^{(m)},\iy)\to\xi_t(x,\iy)$
for every $x\in[x_*,x^*]$. Similarly, by left continuity of
$x\mapsto\xi_t[x,\iy)$, and $x_{(m)}\to x$, $x_{(m)}\le x$,
we obtain $\xi_t(x_{(m)},\iy)\to\xi_t[x,\iy)$.
Let
\[
h_{(m)}(x)=\inf\{\phi'(z):z\in[x_{(m)},x^{(m)}]\},
\qquad
h^{(m)}(x)=\sup\{\phi'(z):z\in[x_{(m)},x^{(m)}]\}.
\]
Summing over $k$ in \eqref{60}, we have for every $m$,
\begin{align}
\label{58}
\liminf_n\gamma^n(t)&\ge \frac{1}{\ell}\int_{[x_*,x^*]}h_{(m)}(x)
\frS(\xi_t(x_{(m)},\iy),\xi_t(x^{(m)},\iy))\xi_t(dx)
\\
\label{59}
\limsup_n\gamma^n(t)
&\le\frac{1}{\ell}\int_{[x_*,x^*]}h^{(m)}(x)
\frS(\xi_t(x_{(m)},\iy),\xi_t(x^{(m)},\iy))\xi_t(dx).
\end{align}
It remains to show the expressions of the r.h.s.\ of
\eqref{58} and \eqref{59} both converge to the r.h.s.\ of \eqref{31} as $m\to\iy$.
By bounded convergence, it suffices
that the integrands in both these integrals converge, for every $x$,
to $\phi'(x)\frS(\xi_t[x,\iy),\xi_t(x,\iy))$. The latter convergence
holds because, by continuity of $\phi'$,
$h^{(m)}(x)\to \phi'(x)$ and $h^{(m)}(x)\to \phi'(x)$, and moreover
$\xi_t(x_{(m)},\iy)\to\xi_t[x,\iy)$ and $\xi_t(x^{(m)},\iy)\to\xi_t(x,\iy)$.
This completes the proof.
\qed

\subsection{PDE satisfied by the limit}\label{sec43}

We take limits along an arbitrary convergent subsequence.
By \eqref{24+}, \eqref{57}, Lemma \ref{lem2} and Lemma \ref{lem6},
we have now established that every limit point $\xi$ satisfies \eqref{70}.

Again, fix a convergent subsequence and denote its limit by $\xi$.
Some notation used below is as follows.
For a right-continuous nondecreasing function $V:\R\to[0,1]$
with $V(0-)=0$, let $V(dx)$ denote the Stieltjes measure induced on $\R_+$.
Let
\[
v(x,t)=\xi_t(x,\iy), \qquad \tilde v(x,t)=1-v(x,t)=\xi_t[0,x].
\]
These are right continuous functions for every $t$.
Thus $\tilde v(dx,t)=\xi_t(dx)$.
Next, for $V$ as above, let
the pure jump part and continuous part be denoted, respectively, by
\[
V^\jmp(x)=\sum_{y\in[0,x]}\Del V(y),
\qquad
V^\cts(x)=V(x)-V^\jmp(x).
\]

\begin{lemma}\label{lem7}
$v$ is a weak solution to \eqref{111+}.
\end{lemma}

\proof
By equations \eqref{24+} and \eqref{57} and Lemmas \ref{lem2} and \ref{lem6},
\[
\lan\phi,\xi_t\ran=\lan\phi,\xi_0\ran+\int_0^t\lan b_1\phi'+a\phi'',\xi_s\ran ds
+\Gam(t).
\]
Write $\Gam(t)$ as $b_0\ell^{-1}\int_0^tA(s)ds$ where
\[
A(t)=\int\phi'(x)\frS(v(x-,t),v(x,t))\tilde v(dx,t).
\]
Let $\tilde v(x,t)=\tilde v^\cts(x,t)+\tilde v^\jmp(x,t)$ denote the decomposition
alluded to above. Let also
$U(x,t)=1-v(x,t)^\ell$.
Then
\[
U(dx,t)=\ell v(x,t)^{\ell-1}\tilde v(dx,t).
\]
Hence
\[
U^\cts(dx,t)=\ell v(x,t)^{\ell-1}\tilde v^\cts(dx,t),
\qquad
U^\jmp(dx,t)=\ell v(x,t)^{\ell-1}\tilde v^\jmp(dx,t).
\]
Now,
\[
A(t)=\int\phi'(x)\ell v(x,t)^{\ell-1}\tilde v^\cts(dx,t)
+\int\phi'(x)\frac{v(x-,t)^\ell-v(x,t)^\ell}{v(x-,t)-v(x,t)}\tilde v^\jmp(dx,t).
\]
The second integral on the right is
\[
\sum_{x\in[0,\iy)}\phi'(x)(v(x-,t)^\ell-v(x,t)^\ell)
=\sum_{x\in[0,\iy)}\phi'(x)(U(x,t)-U(x-,t))=\int\phi'(x)U^\jmp(dx,t).
\]
Therefore, using $\phi'(0)=0$,
\begin{align*}
A(t)&=\int\phi'(x)(U^\cts(dx,t)+U^\jmp(dx,t))=\int\phi'(x)U(dx,t)
\\
&=-\int\phi''(x)U(x,t)dx=\int\phi''(x)v(x,t)^\ell dx.
\end{align*}
Using integration by parts we have
$\lan\phi,\xi_t\ran=\int_{\R_+}(v(x,t)-1)\phi'(x)dx$.
Recalling $\phi'=\tilde\phi$ and combining the above calculations, we obtain
\begin{align*}
\int\tilde\phi(x)(v(x,t)-1)dx &= \int\tilde\phi(x)(v(x,0)-1)dx\\
&\hspace{-5em}
+\int_0^t\int (b_1\tilde\phi'(x)+a\tilde\phi''(x))(v(x,s)-1)dxds
+\frac{b_0}{\ell}\int_0^t\int\tilde\phi'(x)v(x,t)^\ell dxds.
\end{align*}
Thus
\begin{align*}
\int \tilde\phi(x)v(x,t)dx &= \int \tilde\phi(x)v(x,0)dx
+\int_0^t\int (b_1\tilde\phi'(x)+a \tilde\phi''(x))v(x,s)dxds\\
&\qquad
+a\int_0^t\tilde\phi'(0)ds
+\frac{b_0}{\ell}\int_0^t\int \tilde\phi'(x)v(x,t)^\ell dxds.
\end{align*}
According to Definition \ref{222}, this shows that
$v$ is a weak solution of \eqref{111+} once it is verified
that $v\in\Ll^\iy_\loc(\R_+,\Ll^1(\R_+))$.

This property can otherwise be stated as
$\sup_{t\in(0,T]}\int_0^\iy \xi_t[x,\iy)dx<\iy$ a.s., for every $T<\iy$.
Invoking Skorohod's representation theorem we may assume
w.l.o.g.\ that, along the chosen convergent subsequence, one has $\bar\xi^n\to\xi$ a.s.
Arguing by contradiction, let there exist a sequence $t_N\in(0,T]$
and an event of positive $\PP$ measure on which $\int_0^\iy\xi_{t_N}[x,\iy)dx\to\iy$. Now, for each fixed $t\in(0,T]$, by Fatou's lemma,
\begin{align*}
\E\int_0^\iy\xi_t[x,\iy)dx
&\le\liminf_n\E\int_0^\iy\bar\xi^n_t[x,\iy)dx
\\
&\le\liminf_n\frac{1}{n}\sum_{i\in[n]}\E\hat X^n_i(t)
\\
&\le\sup_n\max_{i\in[n]}\E\|\hat X^n_i\|^*_T\le c<\iy,
\end{align*}
where $c=c(T)$ and the last assertion follows from Lemma \ref{lem4}.iii.
This contradicts the assumption, and hence the uniform $\Ll^1$ property follows.
\qed

\subsection{Proof of main results}\label{sec44}

We can now prove Theorems \ref{th1}, \ref{th1+}, \ref{th2}, \ref{th3},
\ref{th4} and Propositions \ref{prop1} and \ref{prop2}.

\noi{\bf Proof of Theorem \ref{th1+}.}
Having proved uniqueness in Lemma \ref{lem9},
we proceed to the remaining assertions.
Existence of a $C^\iy(\R_+\times(0,\iy))$ solution is well known;
see e.g.\ \cite[Lemma 3.2, p.\ 68]{godlewski91}.

In the case $\frf(z)=c_1z-bz^\ell$ and $v_0(x)=\xi_0(x,\iy)$,
we have proved in Lemma \ref{lem7} that the function $v$
induced by any subsequential limit $\xi$ is a weak solution.
By uniqueness it must be equal to the smooth function mentioned
above. Since $x\mapsto v(x,t)$ is automatically nonincreasing,
the function $u:=-v_x$ is nonnegative and smooth.
Using this in \eqref{80}, an integration by parts gives the first 3 parts of
\eqref{14_} as well as that $u$ integrates to $1$ for each $t>0$.

As for the initial condition stated as
$u(\cdot,t)dx\to\xi_0(dx)$ in \eqref{14_}, again using \eqref{80}
with $-v_x=u$ gives $\int\phi'(x)u(x,t)dx=\int \phi'(x)\xi_0(dx)+O(t)$ as $t\to0$,
which proves that the initial condition is satisfied.
This shows that $u$ thus defined forms a classical solution to \eqref{14_}.
\qed

\noi{\bf Proof of Proposition \ref{prop1}.}
The fact that $u_\stat$ given in \eqref{69} is a stationary
solution is verified directly.
As for uniqueness, for any solution $u$, the function
$v=\int_\cdot^\iy u(x)dx$ takes values in $[0,1]$, is monotone, and satisfies
\[
\frf(v)'+av''=0, \qquad v(0)=1,\qquad \lim_{x\to\iy}v(x)=0.
\]
Integrating we obtain
$\frf(v)+av'=c_1$, where $c_1$ is a constant. Since $\frf(0)=0$,
we obtain by the condition at infinity that $\lim_{x\to\iy}v'(x)=c_1/a$,
which can only hold if $c_1=0$. Thus $v$ must satisfy the ODE
\[
v'=-a^{-1}\frf(v), \qquad \text{ on } \R_+,\qquad v(0)=1.
\]
Since $\frf$ is Lipschitz on $[0,1]$, the result follows by ODE uniqueness.
\qed

\noi{\bf Proof of Theorem \ref{th1}.}
The existence of a classical solution has been shown in the proof
of Theorem \ref{th1+} above.
The property $\sup_{t\in(0,T]}\int xu(x,t)dx<\iy$
follows from the fact that $v\in\Ll^\iy_\loc(\R_+;\Ll^1(\R_+))$.
The uniqueness in the class of $C^{2,1}$ functions satisfying
the above uniform integrability follows directly from uniqueness
of weak solutions to \eqref{111}, via integration by parts.
\qed

\noi{\bf Proof of Theorem \ref{th2}.}
The $C$-tightness of $\bar\xi^n$ stated in Lemma \ref{lem3},
the fact stated in Lemma \ref{lem7} that for any limit point $\xi$,
$v(x,t):=\xi_t(x,\iy)$ is a weak solution of \eqref{111+},
and the uniqueness stated in Lemma \ref{lem8} imply that
$\bar\xi^n$ possesses a deterministic weak limit $\xi$,
and that $\xi\in C(\R_+,\calM_1)$. The assertion that $\xi$ is given by
$\xi_t(dx)=u(x,t)dx$, where the latter is the unique classical solution
to \eqref{14_}, follows now from Lemma \ref{lem7} and Theorem \ref{th1+}.
\qed

\noi{\bf Proof of Proposition \ref{prop2}.}
First we show that, under the strengthened moment condition, for $t$ fixed,
\begin{equation}\label{302}
\sup_n\max_{i\in[n]}\E[\hat X^n_i(t)^{2+\eps}]<\iy.
\end{equation}
We do this by reviewing the proof of \eqref{65} in Lemma \ref{lem4}.iii.
Thanks to the assumption that $\Phi_\ser$ has a finite $2+\eps$ moment,
the moment estimate \eqref{303} can now be strengthened to
\[
\sup_n\E[(\|\hat S^{0,n}_1\|^*_t)^{2+\eps}]<\iy,
\]
again by \cite[Appendix 1]{krichagina1992diffusion}. With this and the further
$2+\eps$ moment assumptions on the initial condition, the proof of \eqref{302}
now follows very closely to that of \eqref{65}. We omit the details.

Next, \eqref{302} implies, using Minkowski's inequality,
that $n^{-1}\sum_i\hat X^n_i(t)^2=\int x^2\bar\xi^n_t(dx)$
are uniformly integrable. Hence it suffices to show
$\int x^2\bar\xi^n_t(dx)\to\int x^2\xi_t(dx)$ in probability,
and a similar statement for the empirical first moment. We have from Theorem \ref{th2}
that $\bar\xi^n_t\to\xi_t$ in probability in $\calM_1$.
Using the Skorohod representation theorem, we assume without loss of generality
that this convergence holds a.s. Then $\int x^2\bar\xi^n_t(dx)\to\int x^2\xi_t(dx)$ a.s.\ will
follow once the uniform integrability condition
$\sup_n\int x^{2+\eps}\bar\xi^n_t(dx)<\iy$ a.s.\ is verified.
But this follows from Fatou's lemma in view of \eqref{302}.
\qed

\noi{\bf Proof of Theorem \ref{th3}.}
By Lemma \ref{lem1}.iii and Lemma \ref{lem4}.i, for $i\in[k]$,
\begin{equation}\label{096}
\hat X^n_i(t)=\hat X^n_i(0-)+\hat E^n_i(t)
-\hat S^n_i(T^n_i(t))+\hat b^n_1t
+\hat\la^n_0n^{-1/2}\int_0^tp_{n,\calR^n_i(s)}ds+\hat M^{A,n}_i(t)
+\hat L^n_i(t),
\end{equation}
and $\int\hat X^n_i(t)d\hat L^n_i(t)=0$.
From Step 5 in the proof of Lemma \ref{lem2} we have that
$T^n_i\to\iota$ in probability.
Recall that $(\hat E^n_i,\hat S^n_i)\To(E_i,S_i)$ and
that $E_i-S_i$ are equal in law to $\sig_iW_i$, where $W_i$
are mutually independent standard BMs.
Since $(\hat X^n_i(0-))_{i\in[k]}\To (X_i(0))_{i\in[k]}$ by assumption,
using the dependence structure \eqref{095}, we obtain
\[
(\hat X^n_i(0-)+\hat E^n_i-\hat S^n_i(T^n_i))\To
(X_i(0)+\sig_iW_i)_{i\in[k]},
\]
where $(X_i(0))_{i\in[k]})$ is independent of $(W_i)_{i\in[k]}$.

From \eqref{18} in the proof of Lemma \ref{lem4} we have that
$\hat M^{A,n}_i\to0$ in probability. From Lemma \ref{lem4}.iii,
we have $C$-tightness of $(\hat X^n_i,\hat L^n_i)_{i\in[k]}$.
If we denote the integral term in \eqref{096} by $\hat K^n_i$, we have
tightness of the tuple $(\hat X^n_i,\hat X^n_i(0-)+\hat E^n_i-\hat S^n_i(T^n_i),
\hat K^n_i,\hat L^n_i)_{i\in[k]}$, and denoting a subsequential weak limit point by
$(X_i,X_i(0)+\sig W_i,K_i,L_i)_{i\in[k]}$, one has, for $i\in[k]$,
\[
X_i(t)=X_i(0)+\sig W_i(t)+b_1t+K_i(t)+L_i(t),
\qquad \int X_idL_i=0.
\]
By uniqueness in law of the system of SDE \eqref{100},
the result will be proved once it is shown that
$K_i(t)=b_0\int_0^tv(X_i(s))^{\ell-1}ds$.
Using Skorohod's representation theorem we assume without loss that
the convergence along the subsequence is a.s.
Thus, in view of \eqref{097}, it suffices to show that, along the subsequence,
one has
\begin{equation}\label{099}
(\bar\calR^{c,n}_i)_{i\in[k]}\to(v(X_i(\cdot),\cdot))_{i\in[k]},
\end{equation}
a.s., in the uniform topology on $[t_0,T]$, for any $0<t_0<T$.

To this end, recall that Theorem \ref{th2} establishes, in particular, that
$\bar\xi^n\to\xi$ in $D(\R_+,\calM_1)$, where $\xi\in C(\R_+,\calM_1)$,
and $\xi_t$ is atomless for every $t>0$.
Again it may be assumed that the convergence is a.s.
Hence, with $d_{\rm L}$ the Levy-Prohorov metric,
$\|d_{\rm L}(\bar\xi^n,\xi)\|^*_T\to0$ a.s.
Because for any $t_0>0$ and $x_0>0$,
$v$ is uniformly continuous on $[0,x_0]\times[t_0,T]$, this gives
\begin{equation}\label{098}
\sup_{t_0\le t\le T}\sup_{x\le x_0}|\bar\xi^n_t[x,\iy)-\xi_t[x,\iy)|\to0, \quad \text{a.s.}
\end{equation}
By the atomless property of $\xi_t$, $t>0$, we know that
\[
\bar k_n:=\sup\{\bar\xi^n_t(\{x\}):(x,t)\in[0,x_0]\times[t_0,T]\}\to0,
\quad \text{a.s.}
\]
Now, $\bar\calR^{c,n}_i=n^{-1}(n-\calR^n_i)$ satisfies
\[
\bar\xi^n_t[\hat X^n_i(t),\iy)-\bar k_n
\le\bar\calR^{c,n}_i(t)
\le\bar\xi^n_t[\hat X^n_i(t),\iy)+\bar k_n,
\qquad (x,t)\in[0,x_0]\times[t_0,T],\ i\in[k].
\]
Using this in \eqref{098} shows
\[
\sup_{t_0\le t\le T}|\bar\calR^{c,n}_i(t)-v(\hat X^n_i(t),t)|1_{\{\hat X^n_i(t)\le x_0\}}\to0,
\quad \text{a.s.}
\]
Hence \eqref{099} follows by the a.s.\ convergence $\hat X^n_i\to X^n$
and a.s.\ finiteness of $\|X_i\|^*_T$, $i\in[k]$, completing the proof.
\qed

\noi{\bf Proof of Theorem \ref{th4}.}
The assumed exchangeability of $\hat X^n_i(0-)$ for every $n$, along with
the convergence in probability of their empirical law $\bar\xi^n_{0-}$
to the deterministic limit $\xi_0$, imply the convergence
$(\hat X^n_i(0-))_{i\in[k]}\To(X_i(0))_{i\in[k]}$
where the latter are mutually independent and each is $\xi_0$-distributed
\cite[Proposition 2.2]{sznitman1991topics}. Hence the hypotheses of Theorem \ref{th3}
hold. Because Theorem \ref{th3} asserts that $(\hat X^n_i,\hat L^n_i)_{i\in[k]}$
converge in law to the solution of the system \eqref{100}, in which
$(W_i)$ are mutually independent and are independent of $(X_i(0))$,
the additional mutual independence of $(X_i(0))$ that has just been shown
completes the proof.
\qed

\begin{remark}\label{rem6+} (Sampling with replacement).
If instead of sampling without replacement the algorithm employs
sampling with replacement, one has the expression
\[
\tilde p_{n,r}=\Big(\frac{n-r+1}{n}\Big)^\ell-\Big(\frac{n-r}{n}\Big)^\ell
\]
instead of $p_{n,r}$ from \eqref{19}.
The proofs remain valid because the asymptotic of $\tilde p_{n,r}$
is as that of $p_{n,r}$. More precisely, analogously to \eqref{17},
we have $\max_{r\in[n]}\tilde p_{n,r}=\frac{\ell}{n}+o(\frac{1}{n})$.
Moreover, it is easy to check that \eqref{097} holds true for $\tilde p_{n,r}$.
As a result, the proof of Lemma \ref{lem4}.iii, where \eqref{17} is used,
and of Lemma \ref{lem2} (Step 6), where \eqref{097} is used,
hold verbatim. These are the only places where the expression
for $p_{n,r}$ is used in the proofs.
\end{remark}

\noi{\bf Acknowledgement.}
The first author was supported by ISF (grant 1035/20).

\bibliographystyle{is-abbrv}

\bibliography{main}

\vspace{.5em}

\end{document}